\documentclass[aip, cha, amsmath,amssymb,
floatfix,			
nofootinbib,		
onecolumn,			
preprint,
longbibliography,	
]{revtex4-2}

\usepackage{multirow} 

\usepackage[total={6.5in,8.75in}, top=1.2in, left=0.9in, includefoot]{geometry}
\usepackage{graphicx}
 \usepackage[T1]{fontenc}    

\usepackage[pdftex,bookmarks, colorlinks, citecolor =blue, breaklinks]{hyperref}
\usepackage{url}
\usepackage{xcolor}

\newcommand{\Eq}[1]{(\ref{eq:#1})}

\newcommand{\Lem}[1]{Lem.~\ref{lem:#1}}

\newcommand{\Sec}[1]{\S \ref{sec:#1}}
\newcommand{\Fig}[1]{Fig.~\ref{fig:#1}}
\newcommand{\Figs}[2]{Figs.~\ref{fig:#1}-\ref{fig:#2}}
\newcommand{\Tbl}[1]{Table~\ref{tbl:#1}}
\newcommand{\App}[1]{App.~\ref{app:#1}}

\newcommand{\WB} {\mathit{WB}}
\newcommand{\Case}[1]{\textit{Case}~({#1})}

\newcommand{\InsertFig}[4]
{\begin{figure}[h!t]
       \centerline{
         \includegraphics[width=#4\columnwidth]{./figures/#1}
       }
       \caption{{\footnotesize  #2}
       \label{fig:#3}}
\end{figure}}

\newcommand{\InsertFigTwo}[5] {
\begin{figure}[h!t]
       \centerline{
         \includegraphics[width=#5\columnwidth]{./figures/#1}
         \hskip 0.5in
         \includegraphics[width=#5\columnwidth]{./figures/#2}
       }
       \caption{{\footnotesize  #3}
       \label{fig:#4}}
\end{figure}}

\newcommand{\InsertFigFour}[7] {
\begin{figure}[h!t]
       \centerline{
\renewcommand{\arraystretch}{0.01}
         \begin{tabular}{cc}
         \includegraphics[width=#7\columnwidth]{./figures/#1}&  \includegraphics[width=#7\columnwidth]{./figures/#2} \\
        \includegraphics[width=#7\columnwidth]{./figures/#3}  &  \includegraphics[width=#7\columnwidth]{./figures/#4}
        \end{tabular}
       }
       \caption{{\footnotesize  #5}
       \label{fig:#6}}
\end{figure}}

\newcommand{\bN}{{\mathbb{ N}}}
\newcommand{\bQ}{{\mathbb{ Q}}}
\newcommand{\bR}{{\mathbb{ R}}}

\newcommand{\bT}{{\mathbb{ T}}}
\newcommand{\bZ}{{\mathbb{ Z}}}

\newcommand{\cR}{{\cal R}}

\newcommand{\eps}{\varepsilon}

\newcommand{\digT}{\mathop{\mathrm{dig}}\nolimits_T}
\newcommand{\epsC}{\eps_{\rm crit}}

\newtheorem{thm}{Theorem}
\newtheorem{lem}[thm]{Lemma}

\newcommand{\beq}[1]{\begin{equation}\label{eq:#1}}
\newcommand{\eeq}{\end{equation}}

\newenvironment{se}[1]{\equation\label{eq:#1}\aligned}{\endaligned\endequation}
\newcommand{\bsplit}[1]{\begin{se}{#1}}
\newcommand{\esplit}{\end{se}}

\begin{document}

\title{Proportions of Incommensurate, Resonant, and Chaotic Orbits for Torus Maps}

\author{E. Sander}
\affiliation{Department of Mathematical Sciences,
George Mason University,
Fairfax, VA 22030, USA}
\email{esander@gmu.edu}

 \author{J.D.~Meiss} 
\affiliation{Department of Applied Mathematics,
University of Colorado,
Boulder, CO 80309-0526, USA}
\email{jdm@colorado.edu}

\date{\today}

\begin{abstract}
This paper focuses on distinguishing classes of dynamical behavior for one- and two-dimensional torus maps, in particular between orbits that are incommensurate, resonant, periodic, or chaotic. We first consider Arnold's circle map, for which there is a universal power law for the fraction of nonresonant orbits as a function of the amplitude of the nonlinearity. Our methods give a more precise calculation of the coefficients for this power law. For two-dimensional torus maps, we show that there is no such universal law for any of the classes of orbits. However, we find different categories of maps with qualitatively similar behavior. Our results are obtained using three fast and high precision numerical methods: weighted Birkhoff averages, Farey trees, and resonance orders.
\end{abstract}


\maketitle

\begin{quotation}
We study nonlinear one- and two-dimensional torus maps, starting with Arnold's circle map.\cite{Arnold61b} Jensen and Ecke
and their collaborators \cite{Jensen84, Ecke89} showed that the proportion of parameters for which the map has a dense
orbit on the circle is given by a power law as a function of the amplitude of the nonlinearity,
up to a critical amplitude where the map becomes noninvertible. We compute the 
power law parameters with higher accuracy. We then classify orbit types
for 2D torus maps with typical nonlinearities.  
Grebogi et. al.~\cite{Grebogi83, Grebogi85} considered such maps, but were only able to give relatively imprecise results. With improved numerical techniques, we show that there is no universal power law for the proportion of regular, nonresonant orbits, in contrast to the 1D case. Instead, we find several categories of nonlinearities  for which there are different behaviors. 
\end{quotation}

\section{Introduction}

The distinction between regular and chaotic orbits is fundamental to the study of dynamical systems.
In this paper we use the efficient and accurate method of weighted Birkhoff averages 
(WBA),~\cite{Das16a, Das17, Das19, Sander20, Meiss21} 
to distinguish between these classes of orbits by differences in the rate of convergence of the average.
The WBA can also compute frequency vectors of regular orbits with high accuracy.
We previously used this to distinguish between rational and commensurate rotation
vectors using Farey trees \cite{Sander20} and resonance orders\cite{Meiss21}
for area- and volume-preserving maps.
In the current paper we study maps on the torus, and show that these three methods lead to a
precise and efficient classification of their orbits as chaotic, resonant, or incommensurate.

Other methods for computing rotation numbers include Laskar's frequency analysis technique,\cite{Laskar93a} which uses a Hanning
window to improve Fourier analysis. It is important to note that \citet{Das18b} proved
that the WBA is super-convergent for Diophantine irrationals, while frequency analysis is only quadratically convergent.
Another technique uses Richardson extrapolation; \cite{Seara06, Luque08, Seara09, Luque14, Villanueva22}
while this  appears numerically to be super-convergent, we are not aware of a proof.
Convergence rates have also used by \citet{Rmaileh15} to distinguish between regular and chaotic behavior. 
More generally, \citet{Alseda21} showed that one can compute rotation intervals for circle maps, 
and  rotation sets of torus maps were computed using set based methods by \citet{Polotzek17}
Other approaches include numerical continuation of invariant tori \cite{Sanchez10}
and explicitly computing the conjugacy to rigid rotation. \cite{Haro16, Blessing23}
For a comparison of some of these methods to the WBA, see the discussion in Das et al.~\cite{Das17}


The paper proceeds as follows. In \Sec{Theory} we give an overview of the theoretical background. 
In \Sec{arnold1D}-\ref{sec:torusMaps} we apply these methods to maps $f: \bT^d \to \bT^d$ for $d=1$ and $2$.
We end with conclusions and future plans in \Sec{Conclusions}. 
The numerical methods we use, which have been developed in other papers, are described in
\App{compmeth}. Appendix~\ref{app:EpsCrit} shows how to compute the critical amplitude
for noninvertibility for $d=2$, which we use in \Sec{torusMaps}.
Finally, \App{parameters} lists the parameters that we have used in our numerical simulations. 

\section{Torus Maps and Rotation Vectors}\label{sec:Theory}
We will consider maps $f: \bT^d \to \bT^d$ that are homotopic to the identity.
In general we can assume that $f$ has the form
\beq{torusMap}
	x' = f(x) = x +\Omega + g(x;a)  \mod 1
\eeq
where $\Omega \in \bT^d$, $a$ is a parameter vector, and the nonlinear term $g$ is periodic,
$g(x+m;a) = g(x;a)$ for any $m \in \bZ^d$ (for every parameter $a$). 
We will study several simple examples.
In \Sec{arnold1D}, we consider Arnold's circle map, where $d=1$, and 
\beq{ArnoldMap}
g(x;a) = \frac{a}{2\pi} \sin(2\pi x).
\eeq
In \Sec{torusMaps} we consider the fully 2D case, where 
$d=2$, and $g_1$ and $g_2$ are both sums of sinusoidal functions.
In all cases, 
\[
	g(x;a) = 0 \mbox{ for } a = 0, 
\]
so that the dynamics of \Eq{torusMap} then reduces to a rigid rotation on $\bT^d$, 
\[
	x(t) = f^t(x(0)) = x(0) + t \, \Omega \mod 1, (a = 0),
\]
in which case $\Omega$ becomes  the frequency or rotation vector.

More generally, to determine the rotation vector for an orbit,
we can lift $f$ to $\bR^d$ using the standard projection $\pi: \bR^d \to \bT^d$.
A map $F: \bR^d \to \bR^d$ is then a lift of $f$ if
\[
	\pi \circ F  = f  \circ \pi , 
\]
Here we take the periods of the torus to be one, so that $F(x) \mod 1 = f( x \mod 1)$.
Since $f$ is homotopic to the identity, $F(x+m) = F(x)+m$ for each $m \in \bZ^d$, 
i.e., the map has degree one. Note that any two lifts, say $F_1$ and
$F_2$, differ by an integer vector $F_1(x) = F_2(x) + m$---indeed this
must be true for any fixed $x$, but by continuity the same vector $m$
must work for all $x$.

The orbit of $x \in \bT^d$ has rotation vector $\omega \in \bT^d$ if the limit
\beq{RotationVector}
	\omega(x, f) = \lim_{t \to \infty} \frac{F^t(x) -x}{t} \mod 1
\eeq
exists. This is independent of the choice of lift; however,
it can depend upon the initial point. 
More general versions of rotation vector can be defined,\cite{Misiurewicz89, MacKay94d}
and sometimes computed,\cite{Polotzek17, Alseda21} but we will only compute \Eq{RotationVector}.

For homeomorphisms of the circle ($d=1$), Poincar\'e proved that strict
monotonicity implies that the limit \Eq{RotationVector} exists and is
independent of $x$. For the form \Eq{torusMap}, this occurs when $|g(x;a)| < 1$ 
for all $x \in \bT$. It was shown by Herman that the
resulting rotation number $\omega$ is a nondecreasing function of
$\Omega$. A circle map $f$ that is smooth and strictly monotonic, $f'(x) >
0$, is a diffeomorphism, and for this case Denjoy showed that
$f$ is topologically conjugate to a rigid rotation when $\omega$ is irrational.\cite{Katok99, Guckenheimer02}

Even when $d >1$ the dynamics of \Eq{torusMap} can be conjugate to a rigid rotation only
if it is a homeomorphism.\cite{Grebogi85} Indeed, suppose that there exists a homeomorphism
$\Phi: \bT^d \to \bT^d$ such that 
\beq{conjugacy}
	f \circ \Phi = \Phi \circ R ,
\eeq
where $R(\theta) = \theta + \omega$ is the rigid translation on $\bT^d$.
Since $f = \Phi \circ R \circ \Phi^{-1}$ is a composition of homeomorphisms, it must be one as well.
Moreover, if $f$ and $\Phi$ are diffeomorphisms, then upon differentiation,
\[
	Df(x) = D\Phi( \Phi^{-1}(x)+\omega) D\Phi^{-1}(x);
\]
therefore, $\det(Df) \neq 0$ for all $x$. This gives the necessary condition:
\begin{lem}\label{lem:Conjugacy}
	If the map \Eq{torusMap} is diffeomorphically conjugate to rigid rotation, then its Jacobian is nonsingular.
\end{lem}

\noindent 
The converse is---apparently---not true when $d>1$.\cite{Kim89} Moreover, nonsingularity of $Df$, 
implying that $f$ is a local diffeomorphism, does not imply that $f$ is a global diffeomorphism.\footnote
{This is related to the Jacobian conjecture, which has not even been proven for polynomials.\cite{Bass82}}
\section{Arnold's Circle Map}\label{sec:arnold1D}
In this section, we study the circle map of \citet{Arnold61b}. That is, we consider \Eq{torusMap} for $d=1$
with \Eq{ArnoldMap}. This family of maps has two-parameters, $(\Omega, a)$, and 
without loss of generality we can assume that $\Omega \in [0,1)$ and $a \ge 0$.
The results are depicted in \Figs{ArnoldHist}{probability1d},
and we quickly summarize them here.
Figure~\ref{fig:ArnoldHist} shows histograms of the precision, $\digT$,
see \Eq{digits} in \App{wba}; the peaks near $\digT \gtrsim 14$ correspond
to nearly double precision accuracy for $\omega_T$ \Eq{omegaWBA} using $T = 10^5$ iterates.
The panels of \Fig{arnold1D} categorize the behavior of orbits for $\Omega \in [0,1)$ and 
$a \in [0,2.5]$ with colors corresponding to $\omega_T$.
Parameters with orbits identified as quasiperiodic using \Eq{IrrationalCriterion}
are shown in \Fig{arnold1D}(a) and  those identified as chaotic using \Eq{1DCriterion} in \Fig{arnold1D}(b).
This data is used to produce \Fig{probability1d},
which shows the proportion of chaotic, periodic, and quasiperiodic orbits as a function of $a$. 

The Arnold map is a homeomorphism when $|a| \le 1$, in which case---as noted in \Sec{Theory}---there
are two possibilities: either the map is conjugate to a rigid rotation \Eq{conjugacy} with irrational $\omega$
so that the orbit is quasiperiodic,
or every orbit is asymptotic to a periodic orbit and $\omega \in \bQ$.
In either case the orbits are not chaotic, and the rotation number \Eq{RotationVector} is independent of $x$.
For $|a|>1$ the orbits can be chaotic. We use the WBA, recalled in \App{wba},
to distinguish between regular and chaotic orbits.
In particular, $T$ iterates are used to calculate the approximation $\omega_T$ \Eq{omegaWBA}
of the rotation vector \Eq{RotationVector}. Upon $T$ additional iterates, we compute $\digT$ \Eq{digits},
an estimate of the precision of $\omega_T$. 
If the precision is low, $\digT < D_T$ \Eq{ChaosThreshold}, the orbit is classified as chaotic. 

Here and for the computations discussed below, we chose a grid of $\Omega \in [0,1)$
that is shifted slightly away from rationals to avoid low-order resonances.
Each orbit begins at the same, arbitrary initial point $x(0)
= 0.117789164297101$ and is initially iterated $500$ times to remove
transients. The new initial point $x(500)$ is then iterated with $T = 10^5$ 
to both compute $\digT$ and $\omega_T$. 
A histogram of $\digT$ from the computation of $\omega_T$ is shown in \Fig{ArnoldHist} for three values of $a$.
In each case there are two well separated peaks in the $\digT$ distributions,
giving a sharp distinction between chaotic and regular orbits. 

For example, when $a = 0.8$, where $f$ is a homeomorphism, the average
accuracy is $\langle \digT \rangle= 15.00$, and the
minimum is $\digT = 6.76$. Only $0.11\%$ of the orbits have $\digT < 9$ and are thus incorrectly
classified as chaotic. By contrast $35\%$ of the orbits have $\digT \ge 16$---\Fig{ArnoldHist}
truncates the accuracy at $16$ since the calculations are in double precision. Moreover, this
tallest peak is trimmed in the figure to make the portions of the histogram with smaller values of
$\digT$ more visible. When $a >1$ there is a third small peak in the distributions of \Fig{ArnoldHist}
near $\digT \sim 8.5$; for these orbits the distinction between regular and chaotic is less clear.
The criterion \Eq{1DCriterion} is conservative in the sense that most of these orbits are taken to be chaotic.

\InsertFig{arnoldHistT10e5}
{Histograms of $\digT$ \Eq{digits} for the Arnold circle map
\Eq{torusMap} with \Eq{ArnoldMap}, for $\omega_{10^5}$ using a grid of
$10^4$ values of $\Omega \in [0,1)$. Three histograms are shown, $a =
0.8$ (blue) $1.5$ (orange) and $2.0$ (yellow). If the difference between
the two averages in \Eq{digits} is no more than $10^{-16}$, we set
$\digT = 16$. Each distribution has a peak at $\digT = 16$ (black); these
are truncated in the figure and have heights $0.35$, $0.65$, and $0.68$,
respectively.}
{ArnoldHist}{0.5}

For regular orbits, we use the Farey tree algorithm of \App{FareyTree} to effectively separate rotation
numbers into rationals and irrationals, thereby distinguishing between periodic and
quasiperiodic orbits. This algorithm computes
the minimal denominator \Eq{qmin} and then uses Criterion \Eq{IrrationalCriterion}
to designate $\omega_T$ as \textit{irrational} or not.
We find that when $a = 0.8$, $60.72\%$ of the $10^4$ orbits shown in
\Fig{ArnoldHist} have effectively irrational rotation numbers while
$39.17\%$ are identified as rational (the remaining $0.11\%$ being
omitted since these have $\digT < 9$). 
Results of these computations for additional values of $a$ are
given in \Tbl{ArnoldFractions}.

\begin{table}[htbp]
   \centering
   \begin{tabular}{@{}lccccc @{}} 
    \hline
   $a$    & $\digT\ge16$	&Chaotic & Rational & Irrational	& \Eq{EckeMeasure} \\
	\hline
	0.5  	&0.1986	&	0.0006	&	0.1919 &	0.8075	&  	0.8044 \\
    0.8		&0.3456	&	0.0011	&	0.3917 &	0.6072	&	0.6033 \\
    0.9		&0.4170	&	0.0009  &	0.5087 &	0.4904  &	0.4853 \\
 	0.99	&0.5974	&	0.0022	&	0.7622 &	0.2356	&	0.2355 \\
 	1.0	  	&0.5095 &	0.0026  &	0.8813 &	0.1161	&  	0		\\ 
 	1.01 	&0.5237	&	0.0580  &	0.9397 &	0.0023  &	\\
    1.02	&0.5326	&	0.0709  &	0.9282 &	0.0009  &	\\
    1.5  	&0.6524	&	0.2343	&	0.7657 &	0.0000  &	\\
    2.0 	&0.6863	&  	0.2711	&	0.7289 &	0.0000	& 	\\
   \end{tabular}
   \caption{Fraction of orbits of Arnold's circle map that have $\digT \ge 16$, are chaotic ($\digT < 9$), and
   have rational or irrational rotation number, using $10^4$ values of $\Omega$
   for each $a$ with $T=10^5$ and \Eq{1DTolerance}.
   For $a \le 1$ compare the irrational fraction (column 5) with the power law \Eq{EckeMeasure} (last column).}
   \label{tbl:ArnoldFractions}
\end{table}

It is interesting to note that, for $a = 0.8$, all but $0.9\%$ rotation numbers with the maximum accuracy
$\digT = 16$, are rational: 
it is easier to compute an accurate value for the rotation number if the orbit is periodic.
On the other hand, $15.7\%$ of the rational rotation numbers do have $9 < \digT < 16$.

The classification of orbits with $(\Omega, a) \in [0,1)\times [0,2.5]$ for $2000 \times 2000$ grid
is shown in \Fig{arnold1D}. When $a < 1$, \Fig{arnold1D}(a) is consistent with the nonresonant orbits having nonzero measure, $\mu(a) > 0$.
\citet{FHL20} used computer-assisted KAM methods to give
rigorous upper and lower bounds on $\mu(a)$.
Their methods gave 
\[
	 0.860 748 < \mu(0.25) < 0.914 161.
\]
The lower bound is their rigorous bound, and the upper was obtained by excluding tongues up to period $20$. 
From our computations $\mu(0.25) = 0.9134$.
Their rigorous computations are very time intensive, and it would be impractical to use them to compute anything like the number of parameter values that we have done. Thus our fast but purely numerical method is a complement to the more time consuming rigorous methods. 

Swiateck \cite{Swiateck88} has rigorously shown that the tongues have full
measure for $a=1$, and Khanin \cite{Khanin91} proved that the Hausdorff
dimension of the nonresonant set is then less than one. 
The computations shown in \Tbl{ArnoldFractions} erroneously predict a small but nonzero fraction of quasiperiodic orbits when $a = 1$, which we attribute to the selected cutoffs in Criterion \Eq{IrrationalCriterion}.

\InsertFigTwo{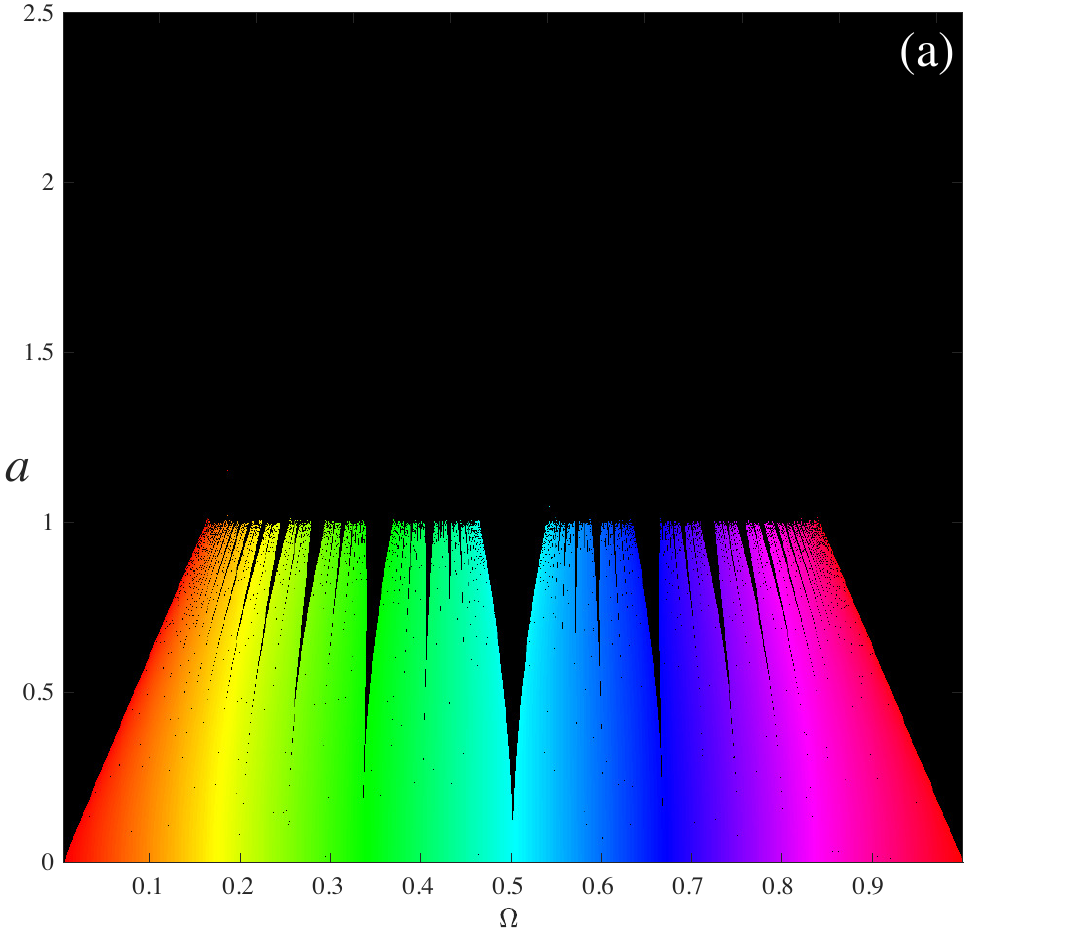}{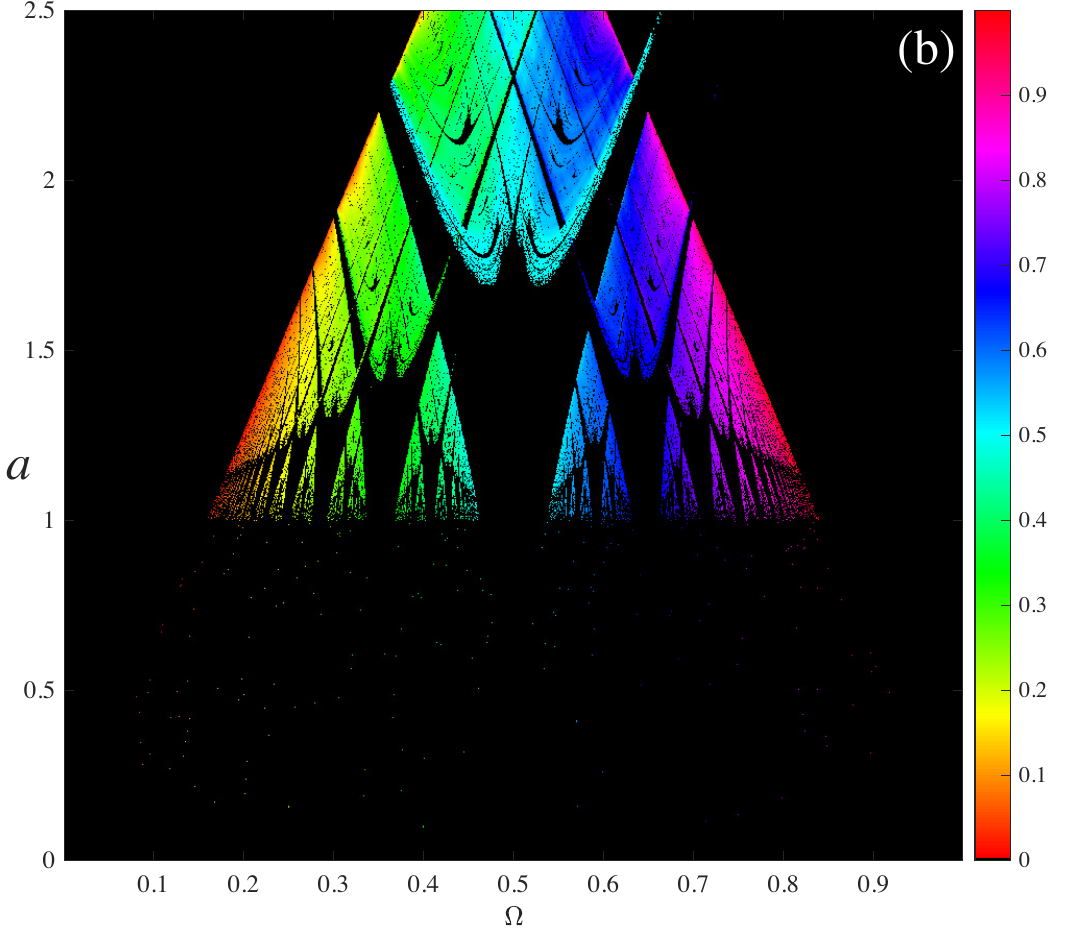}
{Rotation number for the Arnold circle map \Eq{ArnoldMap} as a function
of the parameters $(\Omega,a)$. The color indicates $\omega$, and black indicates no orbits of the given type.
Panel (a) shows the the nonresonant regular orbits and 
panel (b) the chaotic orbits.
}{arnold1D}{0.5}

The computed proportions of chaotic, resonant, and nonresonant orbits are shown 
as a function of $a$ in \Fig{probability1d}. When $0<a<1$, the proportion of nonresonant orbits $\mu(a)$
is close to the previously proposed power law \cite{Jensen84, Ecke89}
\beq{EckeMeasure}
	\mu(a) \simeq (1-a)^{0.314} .
\eeq
As seen in \Tbl{ArnoldFractions}, the value predicted by \Eq{EckeMeasure} is within $0.005$ of our
computations when $a < 1$. For example at $a = 0.99$, the Farey algorithm identifies $2356$ 
rotation numbers as ``irrational'' which is close to $\mu(0.99) = 0.2355$ from \Eq{EckeMeasure}.
We fit the data in \Fig{probability1d} to the more general form
\beq{1DPowerLaw}
	\mu(a) = (1-a)^{p_1  + p_2 (1-a)} 
\eeq
using a log-linear least squares fit.
This form was selected given the known values $\mu(0) =1$, and $\mu(1)=0$, but included a higher order
term in the exponent, since there is no theoretical reason why the power law should have only a single term. 
Our fit gives
\beq{powercoeffs}
	p_1 = 0.3139, \mbox{ and } p_2 = -0.0208.
\eeq 
The root mean squared (rms) error between the
power law and the computed data for $\mu(a)$ (not the log) is 
$0.0024$. By contrast, if we set $p_2 = 0$, the best fit gives $p_1 = 0.3139$ with rms error $0.00449$. 
If we instead compare our data directly with \Eq{EckeMeasure}, the rms error is $0.00454$.

\InsertFig{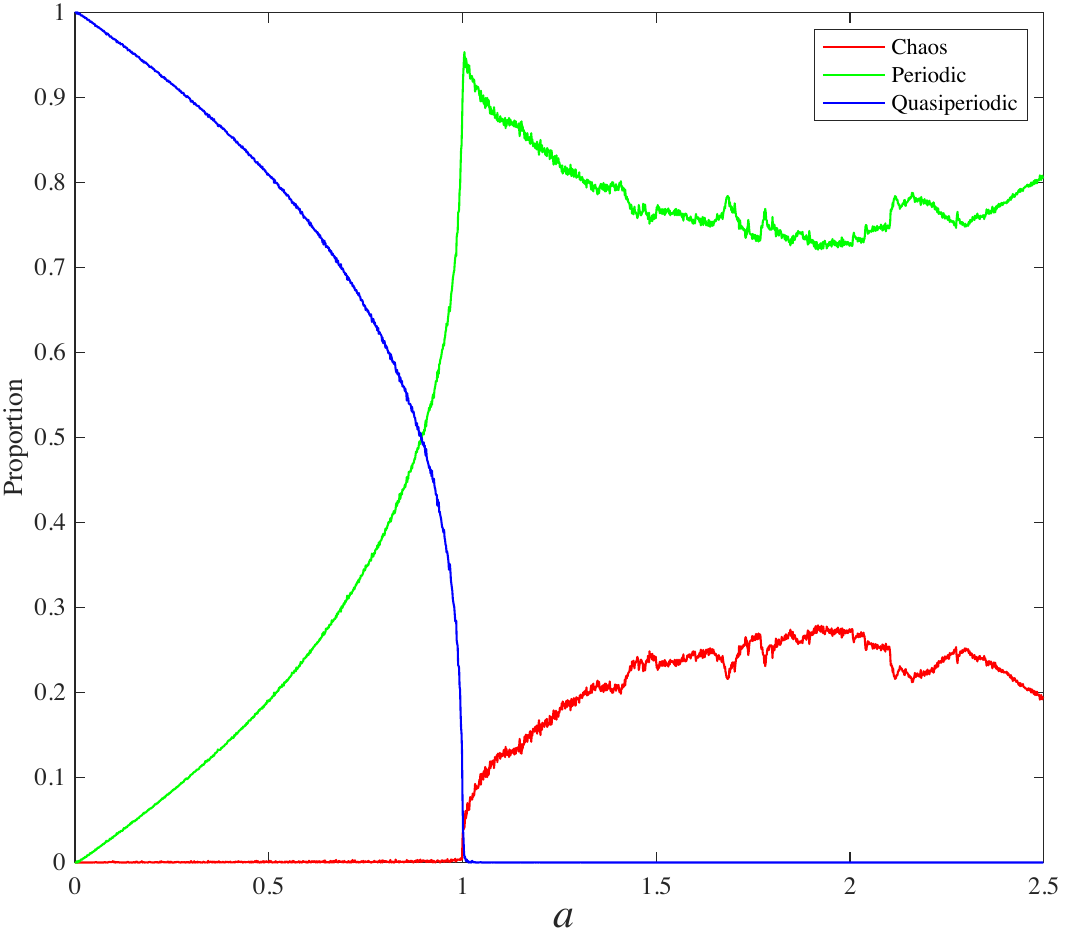}{
The proportion of chaotic, periodic, and quasiperiodic orbits in the Arnold circle map 
as a function of $a$ using the data from \Fig{arnold1D}. For $a \in (0,1)$, the quasiperiodic proportion 
follows the power law \Eq{1DPowerLaw}.}
{probability1d}{0.5}

As is well known, the dynamical behavior of the Arnold circle map changes abruptly at $a=1$.
In \Fig{arnold1D}(b) when $a<1$ there is only a nearly invisible 
``dust'' of points that are falsely labeled as chaotic---$0.07\%$ of the points in this range. 
For $a>1$, the fraction of chaotic orbits in \Fig{probability1d} 
grows but has large fluctuations caused by the well-known bifurcations of the periodic orbits.
Conversely, the computed fraction of quasiperiodic orbits for $a>1$ is essentially  zero.
Indeed when $a=1.02$,  only $9$ of the $10^4$ orbits in in \Tbl{ArnoldFractions} 
are mistakenly identified as ``irrational'' by the Farey algorithm,
and by $a = 1.5$, there are no incorrectly identified orbits. This can also be observed in \Fig{arnold1D}: 
when $a>1$ there are no visible points in panel (a) that would correspond to
falsely labeled nonresonant orbits; moreover,
for $a \in (1.005,2.5]$ the proportion of orbits measured to be nonresonant is less that $0.55\%$.

%
\section{Torus Maps}\label{sec:torusMaps}

We now consider fully coupled maps on $\bT^2$ using the form \Eq{torusMap} with
\beq{forceFunction}
	g(x_1,x_2) = \frac{\eps}{2\pi} \begin{pmatrix}  a_1 \cos(2\pi(x_1 + \phi_1)) +  & a_2 \cos(2\pi(x_2 + \phi_2)) \\
								  a_3 \cos(2\pi(x_1 + \phi_3)) +  & a_4 \cos(2\pi(x_2 + \phi_4)) 
				    		    \end{pmatrix}
\eeq
\citet{Grebogi83, Grebogi85} studied similar maps  in order to gain an understanding of the typical case,
and a number of other specific cases have also been studied.~\cite{Kim89, Ding89, Llibre91, Baesens91, Osinga01, Stark02, Jager06b, Glendinning09}
For simplicity, we will normalize the amplitudes $a \in \bR^4$ so that $\|a \|_1 = 1$; in this case the strength of the forcing function is governed by the parameter $\eps$, and w.l.o.g.~we take $\eps \ge 0$.

In the one-dimensional case there were three types of orbits: periodic ($\omega$ rational),
quasiperiodic ($\omega$ irrational), and chaotic orbits. Since the latter do not occur for a diffeomorphism, verification of the power law \Eq{1DPowerLaw} for Arnold's map when
$0<a<1$ required only the study of a single class of orbit.
In higher dimensions, quasiperiodic orbits can be either resonant or incommensurate,
so there are additional classes of orbits.
As in the circle map case, the proportion of incommensurate orbits is one at $\eps = 0$,
and---as we will see---there exists a critical value $\epsC$ above which there are no incommensurate orbits. 
Thus it seems plausible that a power law like \Eq{PowerLaw} could hold for $\eps < \epsC$.
However since chaotic orbits can occur even for diffeomorphisms, any test of such a form
requires computation of all four classes of dynamical behavior. 

In the first part of this section, \Sec{ResOrdTheory},
we discuss the classification of rotation vectors, and then in
\Sec{EpsCrit} we obtain the critical value $\epsC$ above which the map 
is guaranteed to not be a homeomorphism. The dynamics are studied in \Sec{Typical}
for  a ``typical'' set of amplitudes and phases in \Eq{forceFunction}.
Finally in \Sec{Amplitudes} we study how the proportions of
classes of orbit types vary as the amplitudes change and show that there is no universal
power law.

\subsection{Resonance and Incommensurability}\label{sec:ResOrdTheory}

Perhaps the most natural generalization of the rational verses irrational
dichotomy for rotation numbers to higher dimensions is to ask whether
a vector has rational components, i.e., $\omega = \tfrac{p}{q}$
for some  $p \in \bZ^d$, $q \in \bN$. 
A more general concept, that of \textit{commensurability}, \textit{resonance},
or \textit{mode-locking}, corresponds to the existence of $m \in \bZ^d \setminus\{0\}$ and $n \in \bZ$
such that 
\beq{ResonantPlanes}
	\omega \in \cR_{m,n} = \left\{\alpha \in \bR^d :  m \cdot \alpha =n \right\},
\eeq
a codimension-one plane. Such an $\omega$ has \textit{resonance order} $M = \|m\|_1$ if
this is the smallest length of a (nonzero) vector $m$
for which $\omega \in \cR_{m,n}$.
The set of vectors that do not lie in any resonant plane are \textit{incommensurate}, or {\it nonresonant}.
An example is $\omega = (\sqrt{2},\sqrt{5})$. When $f$ is conjugate
to a rigid rotation \Eq{conjugacy}, these orbits are dense on $\bT^d$. 

For $d=2$, the sets $\cR_{m,n}$ are lines; these are shown in \Fig{ResLines}
up to order $M = 7$. 
The \textit{rank} of a resonance for a given $\omega$ is the number of independent commensurability vectors $m$;
i.e., the dimension of the module of resonance vectors.
In the figure, the rank-two frequency vectors are the points at which nonparallel lines intersect.
Note that $\omega$ is rational only if it has rank two: these correspond to eventually periodic orbits. 

\InsertFig{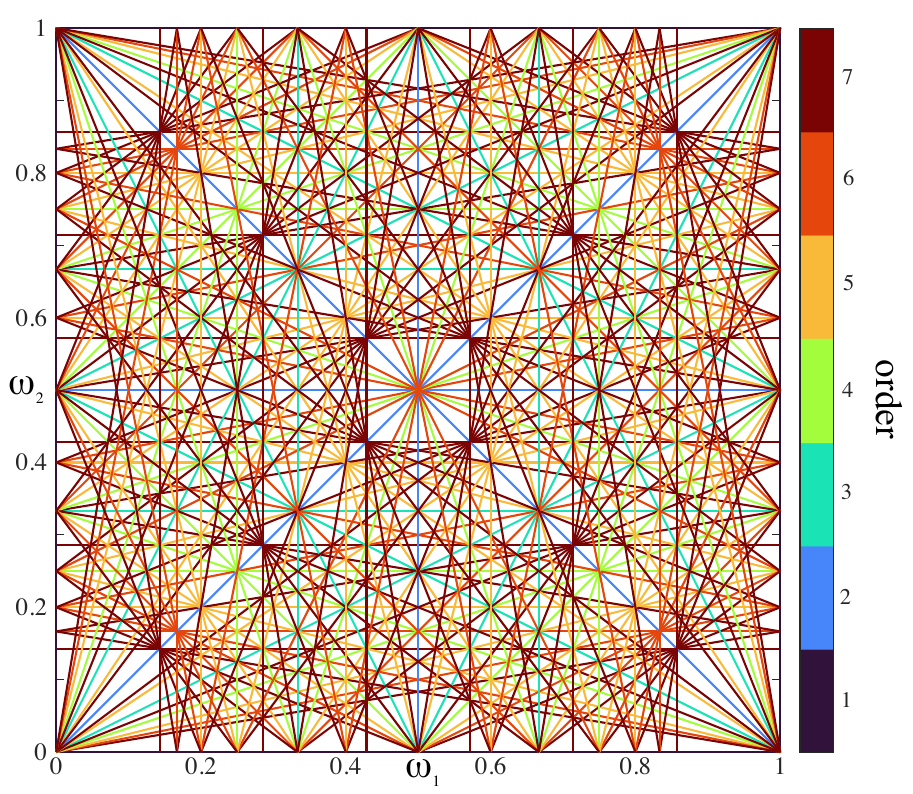}{
Resonance lines up to order $\|m\|_1 = 7$. 
}{ResLines}{0.5}

Commensurabilities that have lower rank are partially
resonant, such as the rank-one vector $\omega = (3\sqrt{2},2\sqrt{2}-1)$ which lies in $\cR_{(2,-3),3}$ so that $M = 5$. A rank-one $\omega$ corresponds to a {\it resonant orbit}, for $d = 2$ these are typically
asymptotic to invariant circles.

The above discussion gives rise to a numerical method for classifying dynamics. 
Each orbit is labeled chaotic or regular using the methods in \App{wba} and Criterion \Eq{2DCriterion}.
Given the accurately approximated $\omega_T$ for a regular orbit, \App{EpsCrit} 
gives a numerical method for distinguishing between the three commensurabilities.
Nonresonant (rank-zero resonant) points, with an incommensurate
rotation vector, are found using Criterion \Eq{Incommensurate}.
Periodic points (rank-two resonant) are found using Criterion \Eq{Periodic}, and 
resonant (rank-one resonant) points are given by to Criterion \Eq{Resonant}.

\subsection{Critical Parameter}\label{sec:EpsCrit}
In this section, we establish the existence of a critical amplitude, $\epsC$, so that 
the map \Eq{torusMap} with \Eq{forceFunction} cannot
be conjugate to a rigid rotation for $\eps > \epsC$: it fails \Lem{Conjugacy}.
The experience with circle maps indicates that local invertibility
should be important dynamically.
Note however that even when $0 \le \eps < \epsC$, where local invertibility holds, 
this lemma does not guarantee the existence of a conjugacy nor the absence of chaotic orbits.

Recall from \Lem{Conjugacy} that a necessary condition for the map \Eq{torusMap} with \Eq{forceFunction}
to be diffeomorphic to rigid rotation is that it has a nonsingular Jacobian,
\bsplit{defineH}
	\det(Df) &= \det \left( I + \eps H  \right) \\
	        &=  \eps^2 \det H +   \eps \; {  \rm tr} H + 1 , \\
    \eps H \equiv Dg &= -\eps \begin{pmatrix}
       a_1 \sin(2 \pi (x_1+\phi_1)) & a_2 \sin(2 \pi (x_2+\phi_2))\\
       a_3 \sin(2 \pi (x_1+\phi_3)) & a_4 \sin(2 \pi (x_2+\phi_4))
    \end{pmatrix}. 
\esplit
To be nonsingular, $\det(Df)$ must be nonzero for all $(x_1,x_2) \in \bT^2$; therefore,
we define the critical value of $\eps$ as
\beq{epsCrit}
	\epsC  \equiv \min \{ \eps \ge 0 : \min_{x \in \bT^2}(\det(Df)) \le 0\} ,
\eeq
i.e., the smallest positive $\eps$ for which $\det(Df) = 0$ for some value of $x$.
When $\eps>\epsC$ conjugacy to rigid rotation is not possible.
Since  $Df = I$ when $\eps = 0$, we know that---if it exists---$\epsC >0$.
In order to compare behavior for different nonlinearities in \Sec{Amplitudes},
we will find that it is appropriate to scale with respect to $\epsC$. 

We claim that if $\det H$ is not identically zero, there is some $x \in \bT^2$ such that 
$\det H <0$. For example, suppose for simplicity that $\phi_i = 0$, $i = 1,\ldots 4$ in \Eq{defineH}.
In this case, since sine is an odd function, if $\det(H(x_1,x_2))>0$, then $\det(H(-x_1,x_2))<0$.
In turn \Eq{defineH} then implies that, if $\eps$ is sufficiently large,
$\det(Df(-x_1,x_2))<0$; therefore, $\epsC$ exists. 
By a similar argument in \App{EpsCrit}, we show that this result holds for 
all phases and amplitudes except for the trivial case $a_1 = a_4 = 0$ and $a_2 a_3 = 0$. 

We find \Eq{epsCrit} numerically using standard root finding methods.
For example, for the coefficients of \Case{0}, given in \Tbl{orbitmparamsC}
of \App{parameters}, we find
\[
	\epsC = 2.22044 .
\]

The existence or nonexistence of a conjugacy to rigid rotation
will not only depend upon $\eps$, but also upon $\Omega$, 
and of course the resulting $\omega$.
Moreover, as we see below, even when $\eps < \epsC$, chaotic orbits can occur.

\subsection{A ``Typical'' Case}\label{sec:Typical}

In this section, we study the dynamics for a fixed set of parameters $a$ 
and $\phi$ as a function of $\eps$ and $\Omega$
(see ``\Case{0}'' in \Tbl{orbitmparamsC} of \App{parameters}). The goal is to
see ``typical'' behavior (as much as one can with a limited set of examples).
While no single parameter set will give all possible dynamics, for \Case{0}
the amplitudes $a_i$ have comparable sizes, and the phases $\phi_i$
are at least not close to rationals with small denominators. We
feel it is useful to look at a single case in detail before giving comparisons between a 
larger set of cases in \Sec{Amplitudes}. 

Six typical phase portraits are shown in \Fig{orbitmontage} for the values of $\eps$ and $\Omega$
listed in \Tbl{orbitmparamsB}.  For first three panels, $\eps = 0.8$. 
Panels (a) and (c) show nonresonant orbits that appear to be dense on $\bT^2$
and for which the rotation vectors, given in \Tbl{orbitmparamsB},
are incommensurate according to Criterion \Eq{Incommensurate}.
Panel (b) shows a resonant orbit
with the low-order resonance $(m,n)=(1,-1,0)$; the orbit lies on an attracting circle that wraps once around
both horizontally and vertically. In panel (d), where $\eps = 1.5$, the orbit is resonant with $m =(2,7,6)$,
as indicated by the fact that the invariant circle wraps seven times horizontally and twice vertically.
Of course, there are also parameters for which the attractor
is periodic for this family of maps, but due to their simplicity, we did not opt to depict any here.
Panels (e) and (f) show examples that Criterion \Eq{2DCriterion} implies are chaotic. Visually (f)
appears to be more chaotic than (e), and a check of the two Lyapunov exponents supports this: for $T = 10^6$,
$\lambda = \{0.0256, -0.0644\}$ for (e) and $\lambda = \{0.2892, -0.0639\}$ for (f).
Though (e) has a positive exponent it seems quite close to the ``weak chaos'' seen for 
quasiperiodically forced circle maps with strange non-chaotic attractors (SNA). We use
the WBA to study such systems in a separate paper.\cite{Meiss24}

\InsertFig{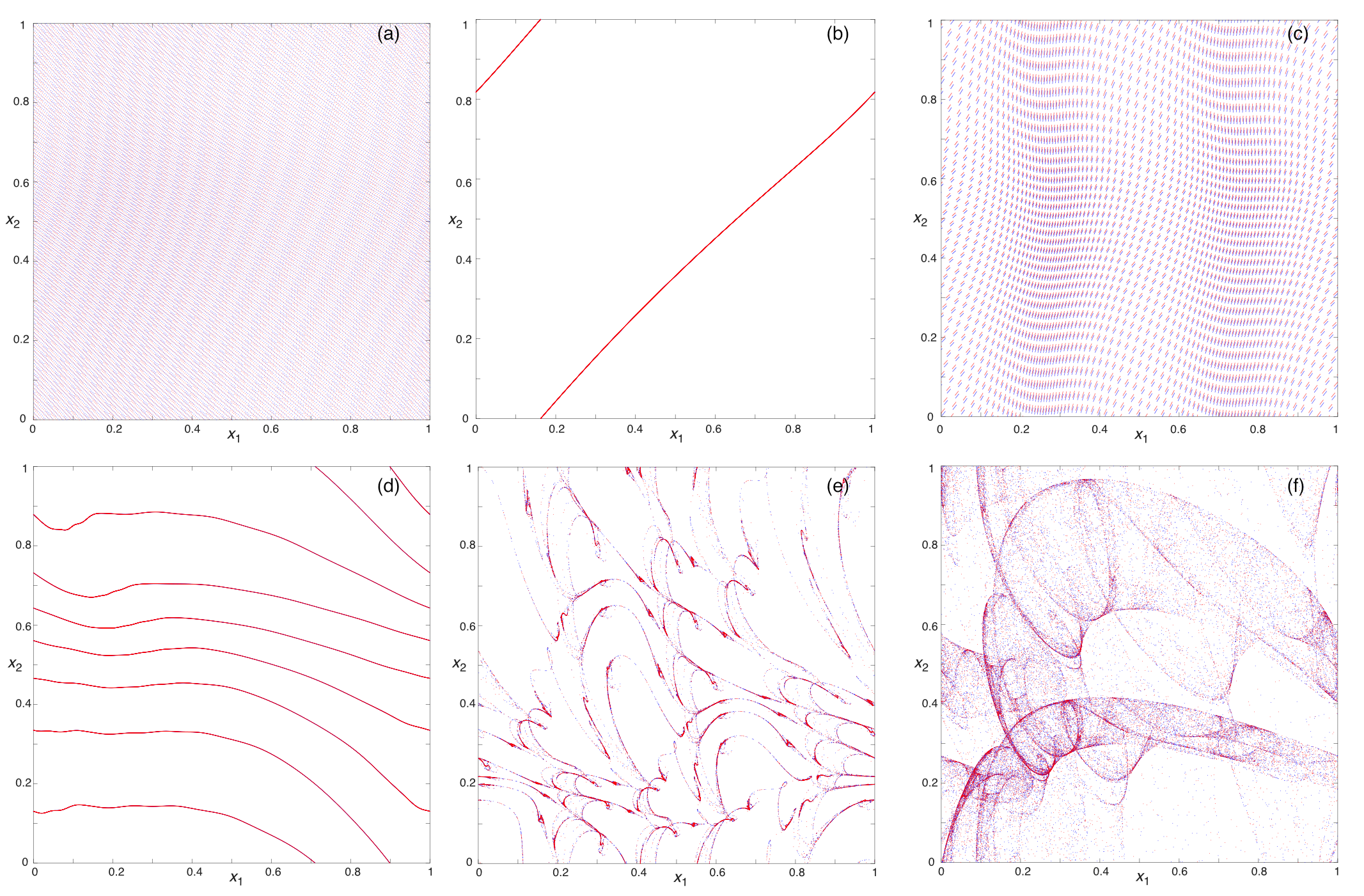}{
Orbits of \Eq{torusMap} on $\bT^2$ with $g$ given by \Eq{forceFunction}. Each panel shows $3(10)^4$ iterates for two 
initial conditions (red and blue), with transients removed.
The amplitudes and phases correspond to \Case{0} in \Tbl{orbitmparamsC},
and values of $\eps$, $\Omega$, the the computed $\omega_T$ and $\digT$ are given in \Tbl{orbitmparamsB}. 
(a) Two-torus;
(b) resonant circle with $(m,n) = (1,-1,0)$;
(c) two-torus;
(d) resonant circle with $(m,n) = (2,7,6)$;
(e,f) chaotic trajectories when the map is noninvertible.
Using $T= 10^6$ for the WBA, gives $\digT > 12$ for the regular orbits in panels (a-d);
for the chaotic orbits (e) $\digT = 4.05$ and (f) $\digT = 2.55$. 
}{orbitmontage}{0.95}

The proportions of the four types of orbits, distinguished using the methods discussed in \App{compmeth},
are shown in \Fig{probability2d} as a function of $\eps$ for \Case{0}.
Here we use a grid of 402 evenly spaced $\eps \in [0,\,1.2\epsC]$
and the same set of 2500 randomly chosen $\Omega \in [0,1)^2$ for each $\eps$. 
Much of the behavior is similar to the 1D case of \Sec{arnold1D}. 
By \Lem{Conjugacy}, the proportion of nonresonant orbits is zero for 
$\eps > \epsC$ (the green dot on the $\eps$-axis in the figure), and indeed we observe that near this point 
the computed proportion (blue curve) does reach zero.
However, unlike the 1D case, this proportion appears to approach zero at $\epsC$ with zero slope.
Moreover, \Fig{probability2d} indicates that chaotic orbits
(red curve) occur for $\eps < \epsC$. 
A similar result, using Lyapunov exponents, was obtained by \citet{Yamagishi20} for high-dimensional torus maps.

We observe that 
resonant proportion (green curve) peaks just below $\epsC$.
The proportion of periodic orbits (purple curve) also 
reaches a maximum but now near $\eps = 2.5$, beyond the resonant peak.
Both the periodic and chaotic proportions grow more-or-less monotonically as $\eps$ crosses $\epsC$.

\InsertFig{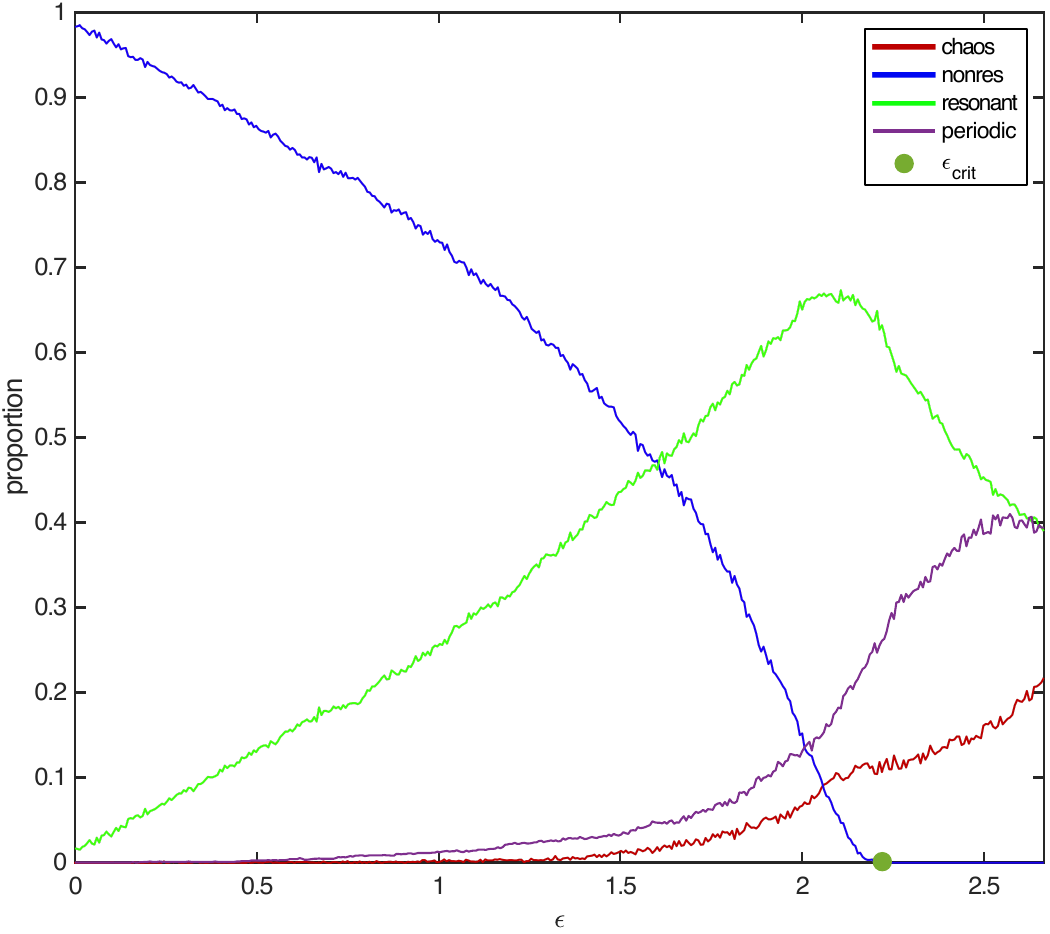}{
The proportion of resonant, nonresonant, periodic, and chaotic orbits in the 2D Arnold map. 
The green dot on the $x$-axis shows the point $\epsC \approx 2.22044$ 
at which the map first becomes locally noninvertible. 
As in the 1D case, there are no nonresonant orbits for $\eps> \epsC$. 
}{probability2d}{0.5}
%

Orbit-type statistics like those shown in \Fig{probability2d} were found in work of
\citet{Grebogi83, Grebogi85} who studied torus maps for $d = 2$ and $3$.
In particular, they considered a finite Fourier series for the
function $g$ in \Eq{torusMap}
and---for randomly chosen amplitudes and phases---computed the fractions of the attractors
that are $m$-tori for $m \in \{0,1,\ldots d\}$, and the fraction that
was chaotic. To do this they computed the Lyapunov spectrum
$(\lambda_1,\lambda_2,\ldots \lambda_d)$. If all $\lambda_i = 0$ the
orbit was classified as lying on a $d$-torus; if all exponents were
nonpositive, and $m$ were zero, the orbit should lie on an $m$-torus.
Finally, if there were any positive exponents, the attractor was classified as chaotic. 
\citet{Yamagishi20} also used the Lyapunov spectrum to classify dynamics of torus maps with $d \sim 100$.
These studies suffer from the problem that accurate computation of the Lyapunov spectrum is difficult.

While \Fig{probability2d} gives a large amount of information on the behavior of orbits aggregated over the 
range $\Omega \in [0,1)^2$, it does not show how these are organized as $\Omega$ varies.
Even when $a$ and $\phi$ are fixed, it would be difficult in a single graph to visualize the categories of orbits
for each $\Omega$ and $\eps$. 
We show slices through this data for fixed $\eps$ in \Fig{arnoldfixepsilon}.
Here the computed $\omega_T$ is shown for nonresonant (top) and resonant and periodic (bottom)
orbits for $\eps = 1.9$ (left), $2.0$ (middle), and $2.1$ (right).
The nonresonant panels show empty strips along resonance lines, recall \Fig{ResLines}, and holes surrounding
the rank-two resonances; these correspond to the Arnold tongues surrounding periodic orbits. 
Kim et al.\cite{Kim85, Kim89} studied the formation of Arnold tongues for a similar torus map using bifurcation theory.
They noted that the region in $(\Omega_1,\Omega_2)$ for which there exists a periodic
attractor (a ``resonance region'') need not be simply connected and that the number of periodic orbits with a given rational
rotation vector $\omega = (p_1,p_2)/q$ can vary. 

\InsertFig{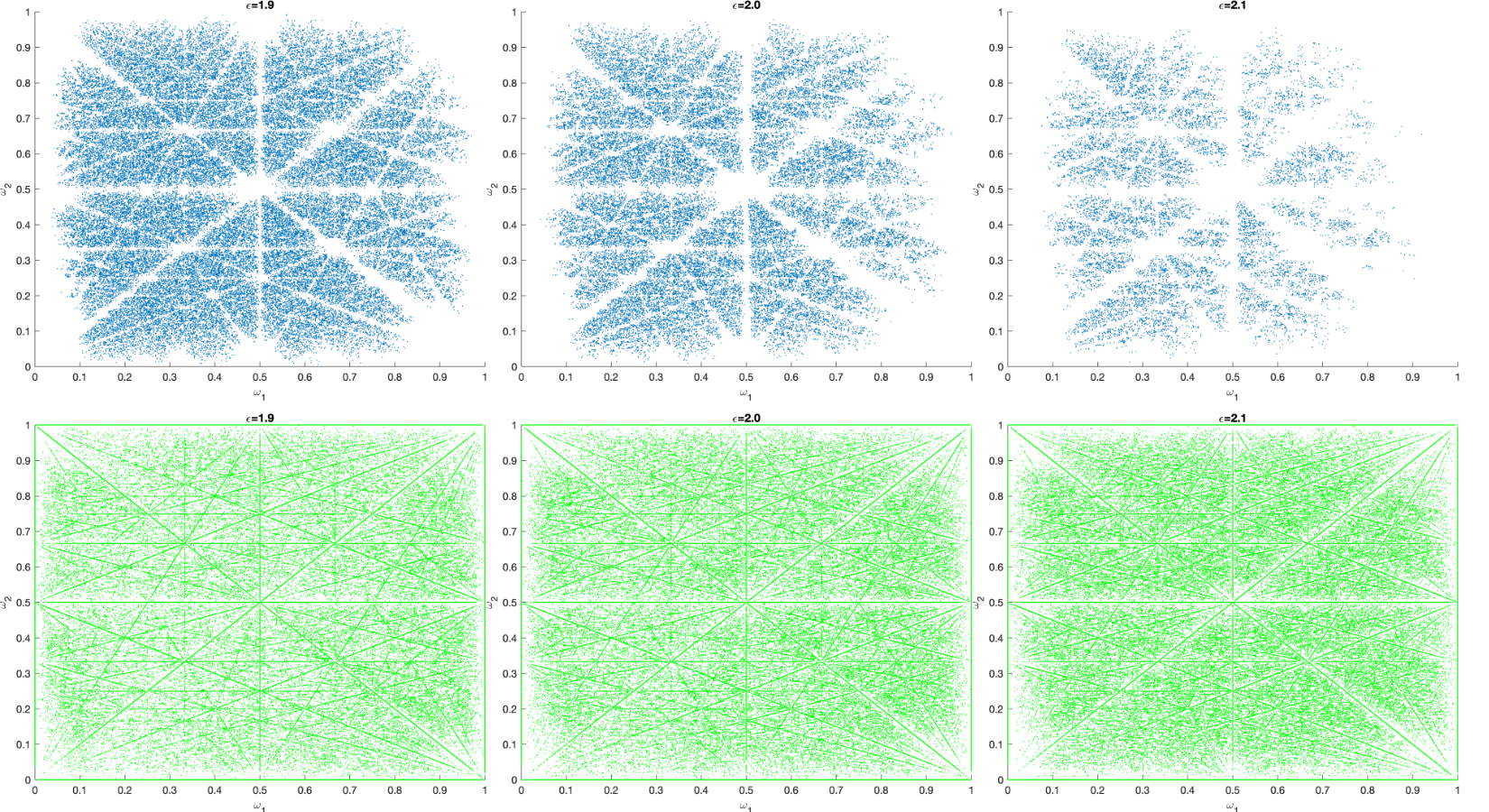}{
The rotation numbers $\omega_T = (\omega_1,\omega_2)$ of the nonresonant (top) and
resonant/periodic (bottom) orbits for \Case{0} in \Tbl{orbitmparamsC} and $250,000$ $\Omega \in [0,1)^2$,
with $\eps = 1.9$ (left), $2.0$ (middle), and $2.1$ (right).
As $\eps$ approaches $\epsC$, the nonresonant set becomes more sparse, and the 
gaps near the lower order resonance lines and periodic points widen.
}{arnoldfixepsilon}{1.0}

To gain further information about the role of $\Omega$, \Fig{singleOmega} shows $\omega_T$ for
resonant (green) and nonresonant (blue) orbits with fixed $\Omega_2 = \gamma = (\sqrt{5}-1)/2$, the golden mean.
This corresponds to a slice through \Fig{arnoldfixepsilon} with fixed $\Omega_2$;
however, now $\eps$ varies over $[0,\epsC]$. This $\Omega$ slice avoids the larger gaps
due to lower-order resonance tongues around the periodic orbits in \Fig{arnoldfixepsilon}. 
Panel (a) is a projection onto the $\omega_T = (\omega_1, \omega_2)$ plane, showing how these vary as $\eps$ grows.
Panel (b) shows the same data, but this time in $(\omega_1,\omega_2,\eps)$; this is
a vector-valued version of a devil's staircase of resonant and incommensurate rotation vectors. 

\InsertFigTwo{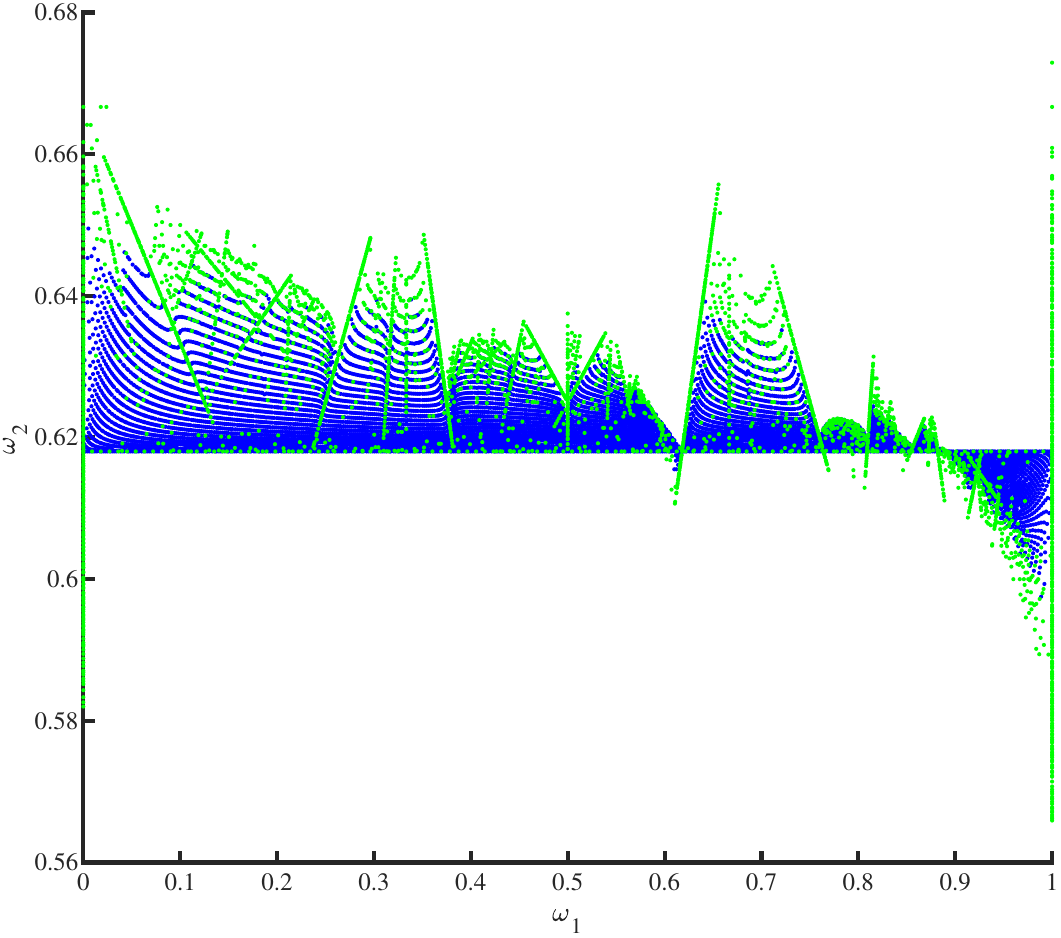}{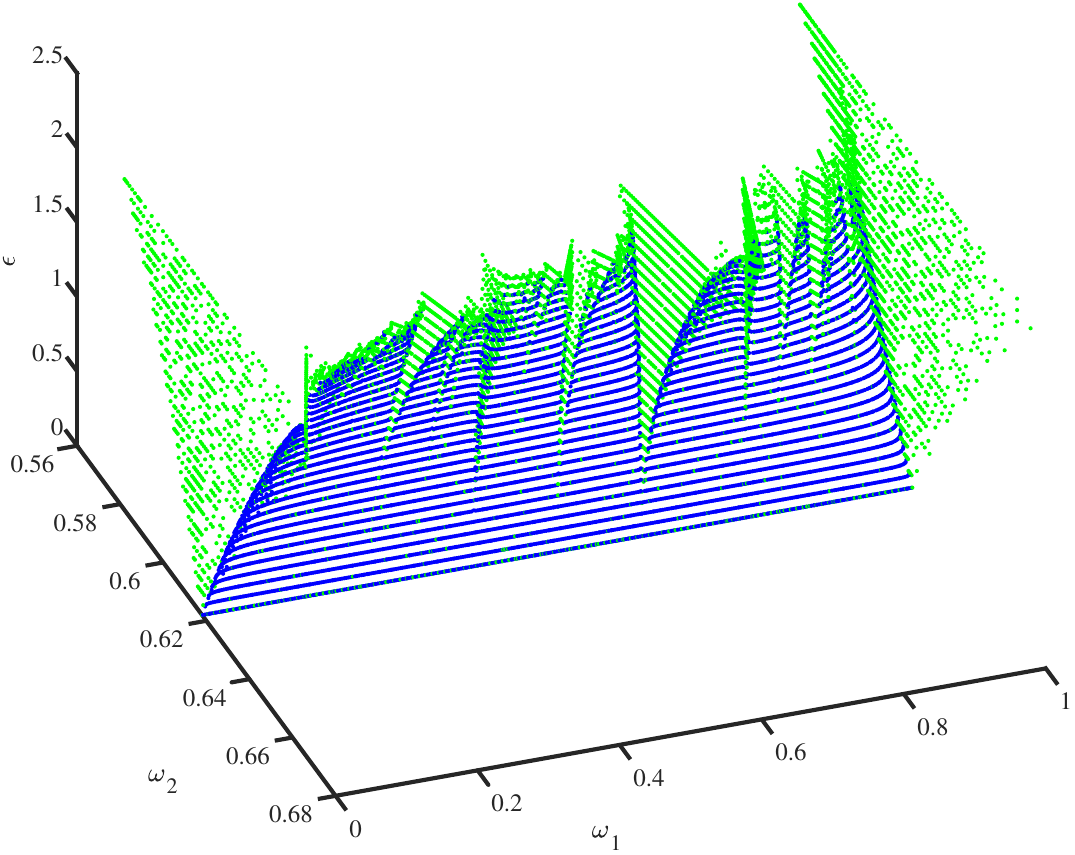}{
(a) Resonant (green) and nonresonant (blue) orbits for $\Omega_1 \in [0,1)$ with $\Omega_2 =\gamma$ 
and 30 evenly spaced values of $\eps \in [0,\epsC]$ showing $\omega_T$ for $T = 10^{6}$.
(b) A 3D view of the data in (a) for coordinates $(\omega_1, \omega_2, \eps)$. 
}{singleOmega}{0.5}

\subsection{Varying the Amplitudes} \label{sec:Amplitudes}

We now proceed to consider more general amplitudes $a$ and phases $\phi$ for \Eq{forceFunction} using
eight parameter sets given in \Tbl{orbitmparamsC} of \App{parameters}. 
The first four sets are chosen randomly,
but the last four are chosen to illustrate noteworthy dynamical categories.
In each case we normalize $\|a\|_1 = 1$ and compute $\epsC$ using \Eq{epsCrit}---this is given in the last column in \Tbl{orbitmparamsC}.
The four panels of \Fig{manycases} show the proportion of (a) nonresonant, (b) resonant, 
(c) periodic and (d) chaotic orbits for the eight parameter sets as a function of $\eps$, scaled by $\epsC$.
One can see immediately from \Fig{manycases} that \Case{4} (purple) and \Case{7} (orange) are outliers.
These correspond to uncoupled and semidirect cases, and will be discussed further below.

Panel (a) shows that in each of the eight cases
the proportion of nonresonant orbits drops to zero at $\eps/\epsC = 1$. 
These curves illustrate that there is no universal power law for the nonresonant orbits of the form
\beq{PowerLaw}
	\mu(\eps) \simeq \left(1 - \frac{\eps}{\epsC} \right)^p \;,
\eeq
that would be analogous to \Eq{EckeMeasure}, since the curves have different shapes. 
Indeed, since some of the curves (e.g., \textit{Cases} (3), (5) and (6)) 
appear to have zero slope at both $\eps = 0$ and $\epsC$,
they cannot not satisfy a single power law on the full range $[0,\epsC]$.
Even if the form \Eq{PowerLaw} was valid asymptotically close to $\epsC$,
$p$ would need to be greater than one for the fully coupled cases,
in stark contrast to \Eq{EckeMeasure}.

Most of the resonant proportions, shown in \Fig{manycases}(b), exhibit smooth peaks for some $\eps < \epsC$.
The outliers are again \textit{Cases} (4) and (7), where there are sharp peaks at $\epsC$.
The periodic orbit proportions peak as well, but not at the same $\eps$ as the resonant peaks
(and in \Case{7} there are no periodic orbits). 
Panel (d) shows that the onset of chaos is considerably below $\epsC$, again with the exception of \textit{Cases} (4) and (7).
All of these observations are consistent with those for other randomly chosen $a$ and $\phi$ that we considered,
but that are not shown.

\InsertFigFour{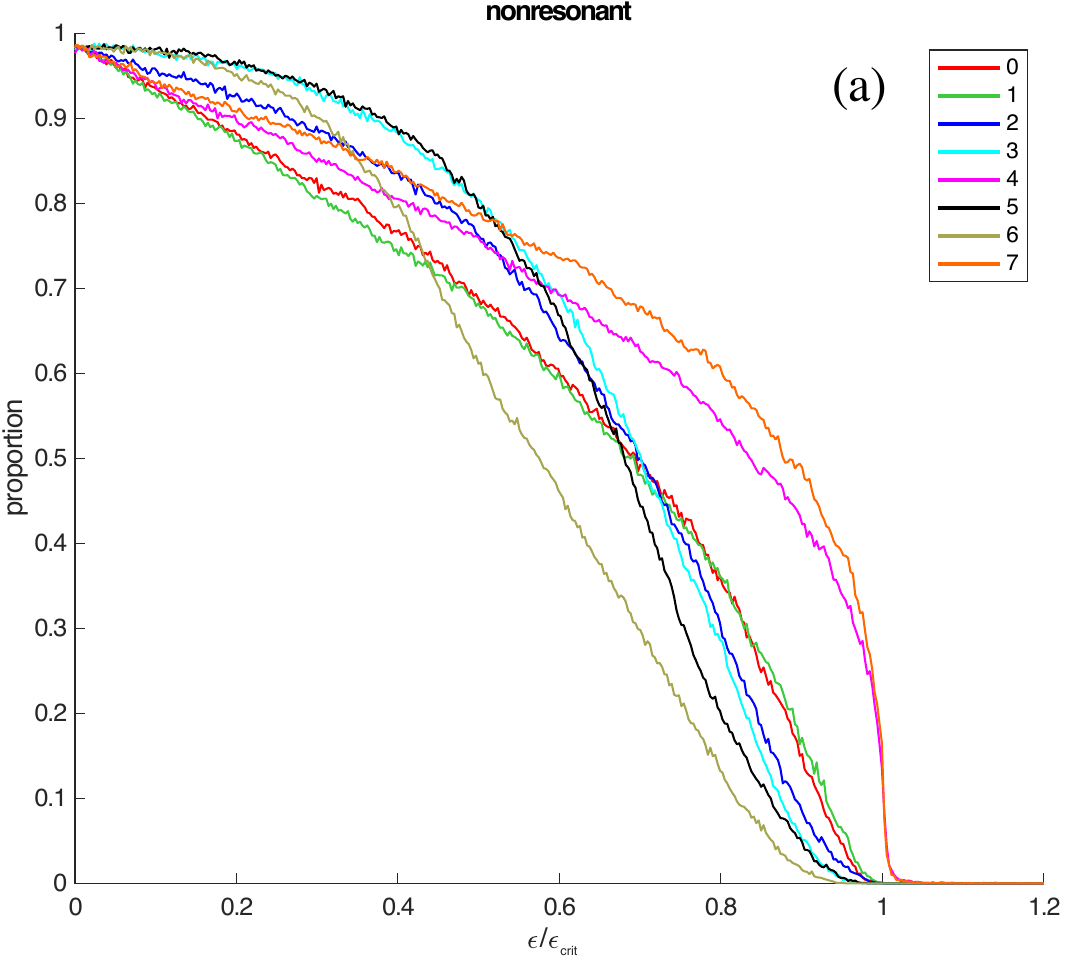}{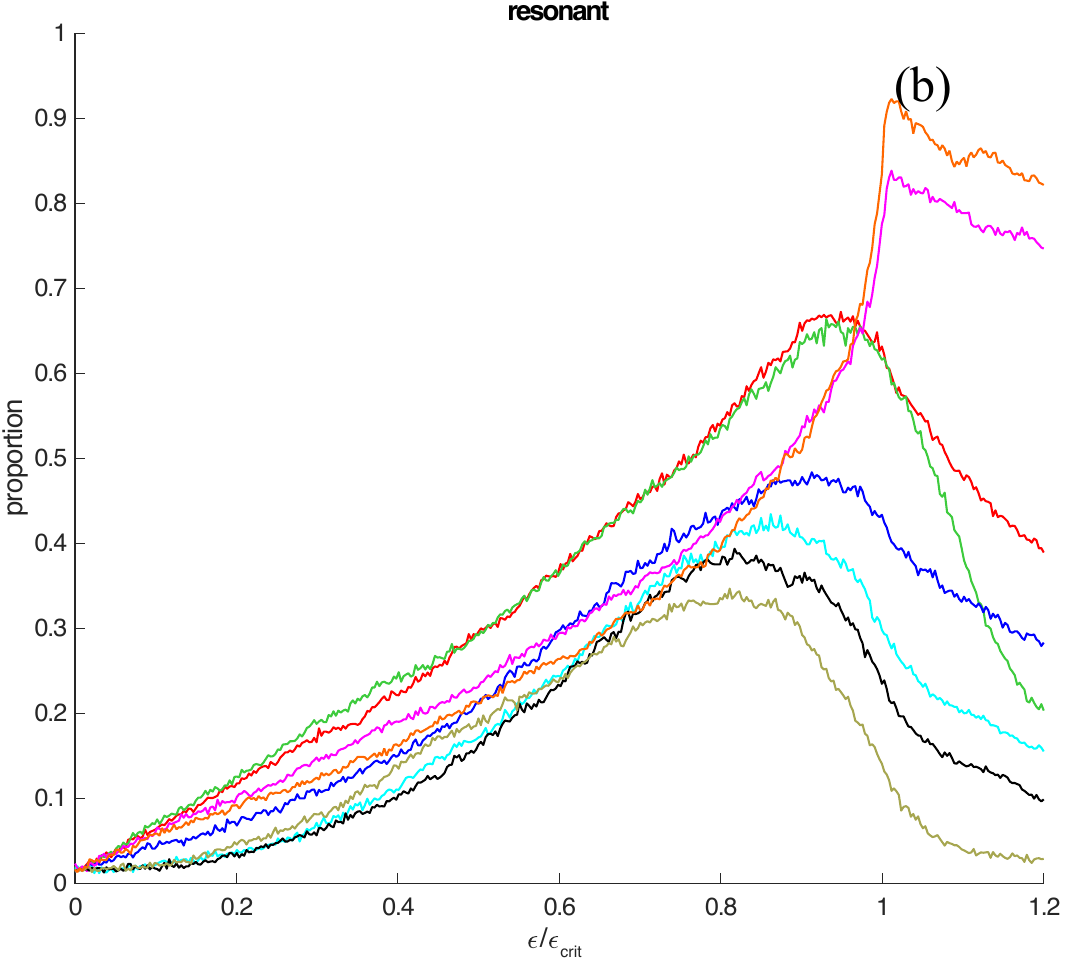}{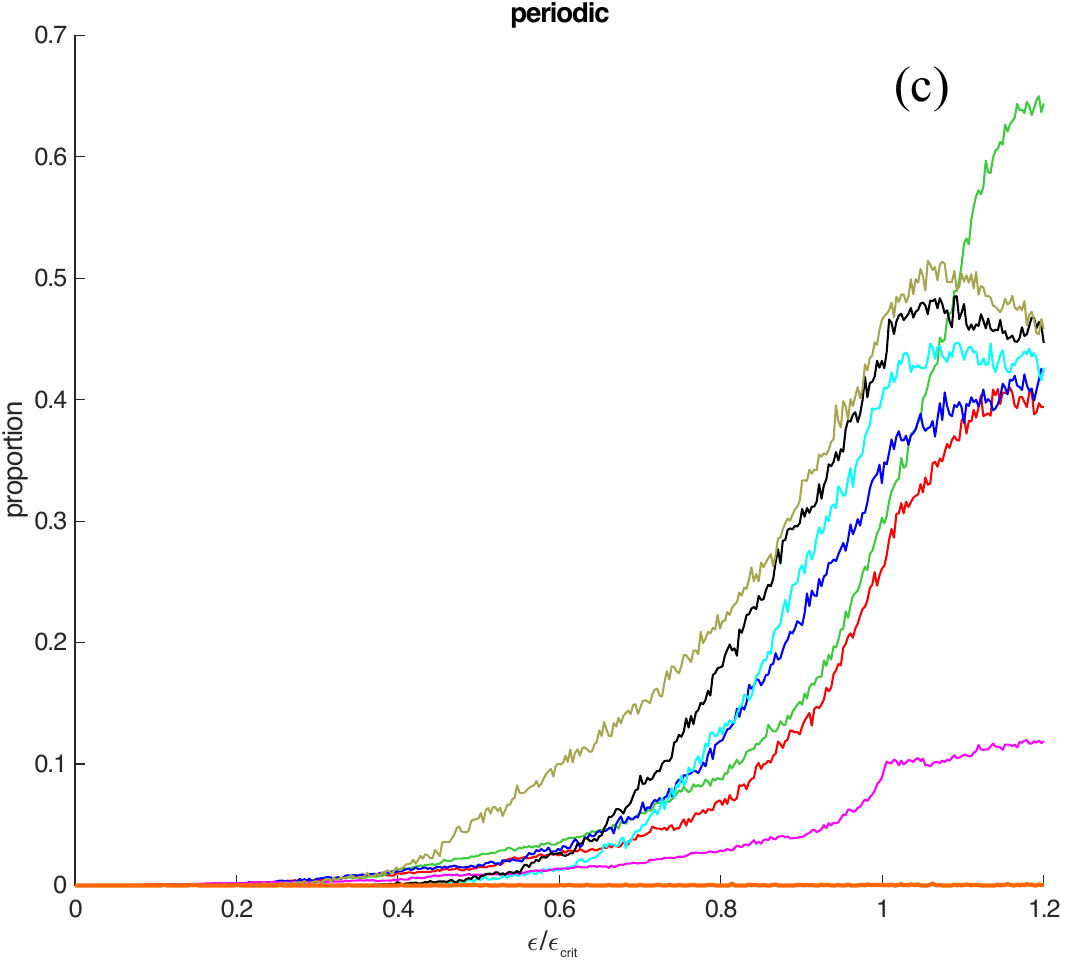}{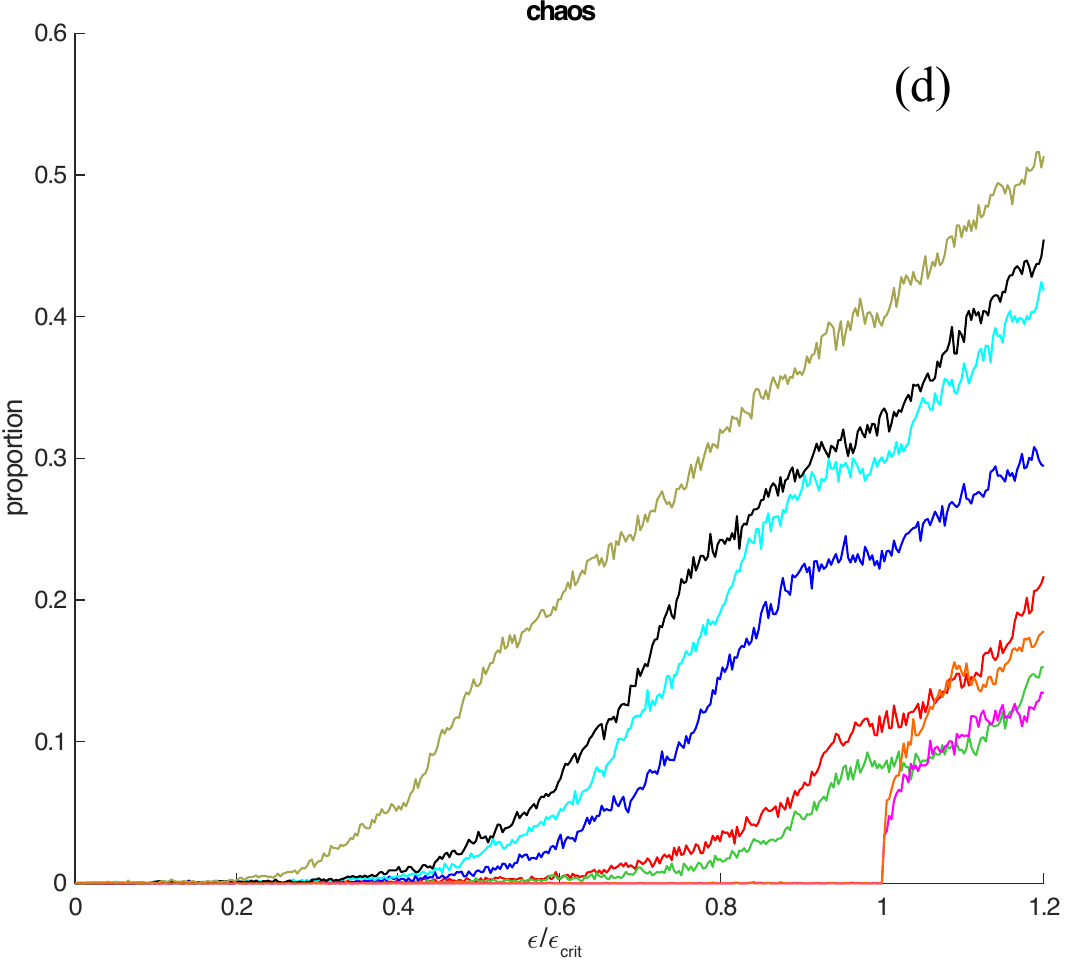}{
The proportion of (a) nonresonant, (b) resonant, (c) periodic, and (d) chaotic orbits as a function of $\eps/\epsC$
for the eight sets of the coupling and phase parameters $a$ and $\phi$ in \Tbl{orbitmparamsC}.
Note that the vertical scales are different for each panel.
}{manycases}{0.4}

We now consider four special subcategories, describing their distinct behaviors.
Though these are extreme examples, we have observed that for nearby amplitudes, the behavior is similar.
Though we have not systematically varied the phases, these do not seem to change the qualitative behavior.

\begin{itemize}
\item {\em Uncoupled Components}: If $a = (1,0,0,0)$ 
the system is uncoupled and essentially 1D, the second  
component $\omega_2 = \Omega_2$ will be irrational for almost all $\Omega$.
Thus the rotation vector is incommensurate depending only the first component which is the Arnold circle map.
Therefore the power law \Eq{PowerLaw} will hold with the 1D value $p_1 \approx 0.3139$ from \Eq{1DPowerLaw}.
Of course the same considerations apply to the case $a = (0,0,0,1)$.

More generally, the uncoupled case corresponds to $a = (1-u,0,0,u)$, and we can restrict to $u \in [0, \tfrac12]$ by symmetry.
Since $u \le \tfrac12$, $\epsC = (1-u)^{-1}$ where the first component becomes noninvertible.
An example is \Case{4} where $u \approx 0.24$. Figure~\ref{fig:manycases} shows that orbit type proportions
for this case have shapes similar to those for Arnold's map in \Fig{probability1d}. 
In particular there is no chaos for $\eps< \epsC$.
Note that the rank-two resonances, panel (c), correspond to the cross product of Arnold tongues
for each of the decoupled 1D maps; there is no analogue of this in \Fig{probability1d}.

For this uncoupled case, the proportion of orbits with irrational $\omega_1$ and that with irrational $\omega_2$ 
will both satisfy \Eq{PowerLaw} but the effective critical parameters will be different, 
$\epsC = (1-u)^{-1}$ and $u^{-1}$, respectively.
However, the incommensurate proportion for the 2D rotation vector $\omega$ can have a different behavior.
The simplest case, $u = \tfrac12$, has $\epsC = 2$ and
would have a nonresonant proportion $\mu$ that is approximately the square of the proportions for each 1D map;
therefore we expect \Eq{PowerLaw} to hold with $p \approx 2p_1$.
However note that this does not take into account rational relations like \Eq{ResonantPlanes}
with both $m_1$ and $m_2 \neq 0$.
By contrast if $u$ is small, the second map will have an irrational proportion that is near $1$ and
changes only slightly as $\eps$ grows to $\epsC = (1-u)^{-1}$.
This should give \Eq{PowerLaw} with $p \approx p_1$, close to the 1D case.
More generally, we observe that the power $p$ varies continuously
between $p_1$ and $2 p_1$ as $u$ grows from $0$ to $\tfrac12$.

We have also looked at weakly coupled examples (not shown); these show orbit-type proportions
similar to \Case{4}.

\item{\em Semidirect product}: When $a = (0,u,0,1-u)$, the behavior is similar to the 1D case
even for nonzero $u$ since the map is a semidirect product. (Of course, the case $(1-u,0,u,0)$ is of this type as well.)
An example is \Case{7} where $a \approx (0,0.35,0,0.65)$.
Here the second component is Arnold's circle map and the first is linear but driven by the dynamics of $x_2$.
Thus this case is also essentially 1D and $\epsC = (1-u)^{-1}$.
Since $a_1 = 0$, we might expect that $\omega_1 = \Omega_1$, 
at least for $\eps < \epsC$ where the dynamics of $x_2$ is either conjugate to an irrational
rotation or is asymptotically periodic.
Note that there are no periodic orbits in panel (c) for this case because the grid of $\Omega$ values
was chosen to avoid the rationals. 

\Case{6} is a weakly coupled perturbation of the trivial semidirect case with $u \approx 1$.
This case is also near the anti-coupled case discussed below.

\item{\em Quasiperiodic forcing}: When $a = (1-u,u,0,0)$ the first component is generically quasiperiodically forced.
We consider this class of maps separately in a forthcoming paper.~\cite{Meiss24}
For fixed $\Omega_2 = \gamma$ the proportion of nonresonant orbits appears to have zero slope at $\epsC$,
and thus does not satisfy a power law with $0<p<1$. 
This also persists when we allowed $\Omega_2$ to vary. Equivalent behavior occurs for $(0,0,u,1-u)$. 

\item{\em Anti-coupled system}: When $a = (0,1-u,u,0)$ the dynamics is fully 2D, but
could be thought of as anti-coupled. For this case $\epsC =(u(1-u))^{-1/2}$.
An example is \Case{5} (black) in \Fig{manycases} where $u \approx 0.23$.
Here we again find that the proportion of nonresonant orbits has zero slope $\eps = \epsC$, with the 
fastest decay when $u \approx 0$ or (by symmetry) $u \approx 1$. 
We suspect that this is due to the larger fraction of chaotic orbits at smaller $\eps/\epsC$. 
Note that \Case{3} is a slightly perturbed anti-coupled system, and its proportion curves are close to \Case{5} even 
though the effective value of $u$ is different.
As mentioned above, \Case{6} can also be thought of as a perturbation of a trivial version of this case with $u=0$.
Mode locking for an anti-coupled case with $u = \tfrac12$ has been previously studied by \citet{Baesens91}.
\end{itemize}

\section{Conclusions and Future Directions}\label{sec:Conclusions}

In this paper, we have used efficient techniques to characterize the dynamics of orbits of
one- and two-dimensional torus maps. 
We used the convergence rate of the WBA to distinguish between chaotic and regular
orbits, defining a threshold for the precision, $\digT$ \Eq{digits}, after a fixed number of iterates $T$.
Given an accurate value for rotation vector $\omega_T$, we determined if the vector is nearly resonant by finding the smallest order of a resonance vector within a distance $\delta$. This allows use
to characterize regular orbits as incommensurate (nonresonant), rank-one (resonant), or rank-two (periodic).
This computation is especially efficient for the 1D case,
where we can use the Farey tree to compute the minimum denominator \Eq{qmin}.\cite{Sander20}
We hope in the future that an efficient algorithm to compute \Eq{Mres} for the 2D case can be found. 
Meanwhile, we used a brute force method, following our previous work on volume-preserving maps.\cite{Meiss21}


Our methods naturally extend to higher dimensional maps.
However, there are some significant computational challenges for $d > 2$. 
The WBA method extends without any difficulty to three and higher dimensions;\cite{Meiss21}
other than the standard problem of needing more initial conditions and iterates to get a sense of the full dynamics,
there are no difficulties to distinguish between chaotic and regular orbits.
Finding periodic (ie. rank-$d$ resonant) orbits in dimension $d$ is also just a matter of
checking that the rotation vector is rational in each component. 

However, there are computational issues that may limit our ability to distinguish lower rank resonances. 
In \App{compmeth}, we observed that for random values in $\bT^1$ the typical denominator of $\delta$-close rational approximation to an irrational scales as $q_{min} \sim \delta^{-1/2}$, recall \Eq{MeanDenom}.
Of course $q_{min}$ is the 1D version of the resonance order $M$, \Eq{Mres}.
In dimension two, for random values in $\bT^2$, we observed that $M$ scales as
$M \sim \delta^{-1/3}$, recall \Eq{MeanRes}.
There is theoretical support for the $d=1$ result \cite{Chen23} and \citet{Marklof24a, Marklof24b} has shown that 
for the general $d$-dimensional case \Eq{MeanRes} becomes
\beq{conjecture}
 	\langle \log_{10}(M(\omega,\delta))\rangle \approx -\tfrac{1}{d+1} \log_{10}\delta + c_d \;.
\eeq

This is a new torment to add to the usual ``curse of dimensionality.'' 
As $d$ grows, the typical resonance order of a theoretically incommensurate vector will be pushed closer to zero
for a given precision $\delta$.
For example, if for $d = 4$ we were to use the same precision, $\delta = 10^{-9}$,
that we used in the current paper, then typical resonance order would be $M \sim 10^{9/5} \sim 63$.
Therefore the computed order of an incommensurate $\omega$ would
be small enough that it would be quite hard to distinguish it from a vector that is actually resonant
with even a modest resonance order.
In order to be able to distinguish these values, one would need the added computational expense of extended precision;
moreover, the computation of a frequency $\omega$ to higher accuracy would also require increasing the
number of iterates $T$. The calculations would quickly become very slow.
We hope that further research in this area will prove us to be overly pessimistic.

\begin{acknowledgments}
ES was supported in part by the Simons Foundation under Award 636383.
JDM was supported in part by the Simons Foundation under Award 601972.
        Useful conversations with Nathan Duignan are gratefully acknowledged.
\end{acknowledgments}

\newpage
\appendix
\begin{center}
\Large{\textbf{Appendices}}
\end{center}


\section{Computational Methods}\label{app:compmeth}

In this section we describe the numerical methods used in \Sec{arnold1D}-\ref{sec:torusMaps}.
We use the weighted Birkhoff average, \App{wba}, for computing rotation number and identifying chaos.
The Farey tree method in \App{FareyTree} distinguishes regular orbits that are periodic
from those that are quasiperiodic for $d=1$. 
The method of resonance orders in \App{ResonanceOrder} extends this to higher dimensions
to distinguish the resonant (lower-dimensional) invariant tori from those that are nonresonant
(full-dimensional).

\subsection{Regularity vs. Chaos: Weighted Birkhoff Averages}\label{app:wba}
We briefly review here the weighted Birkhoff average \cite{Das16a, Das17, Das19} and how it
distinguishes between regular and chaotic orbits.\cite{Sander20} Given a map $f:M \to M$, recall that the time average of a function $h: M \to \bR$ along an orbit of $f$ is simply
\beq{Birkhoff}
	B(h)(z) = \lim_{T \to \infty} \frac{1}{T} \sum_{t=0}^{T-1} h \circ f^t(z) ,
\eeq
if this limit exists. The classic theorem of Birkhoff implies that if the orbit of $f$ is ergodic on a
set with invariant measure $\mu$, and if $h \in L^1(M,\bR)$, then
\[
	B(h)(z) = \langle h \rangle = \int_M h \, d\mu
\]
for $\mu$-almost every $z$. However, the convergence to this limit is at best as $1/T$
and can be arbitrarily slow.\cite{Kachurovskii96, Krengel78}.

To compute the average efficiently and accurately for a length-$T$ segment of an orbit,
we modify \Eq{Birkhoff} using the $C^\infty$ weight function 
\[
	\Psi(s) \equiv \left\{ \begin{array}{ll}  e^{-[s(1-s)]^{-1}}  & s \in (0,1) \\
	         								0	& s \le 0 \mbox{ or } s \ge 1 
						\end{array} \right. \;.
\]
This exponential bump function limits to zero with infinite smoothness 
at $0$ and $1$, i.e., $\Psi^{(k)}(0) = \Psi^{(k)}(1) = 0$ for all $k \in \bN$. 
The finite-time weighted Birkhoff average (WBA) is then defined by
\beq{WB}
	\WB_{T}(h)(z) = \frac{1}{S}\sum_{t = 0}^{T-1} \Psi \left(\tfrac{t}{T}\right) h \circ f^t(z) \;, 
\eeq
with the normalization constant
\beq{SmoothedAve}
	S \equiv \sum_{t=0}^{T-1} \Psi \left(\tfrac{t}{T}\right)  \;. 
\eeq

As shown by \citet{Das16a}, this gives the same answer as $T \to \infty$
as \Eq{Birkhoff}; however, for regular orbits \Eq{WB} can converge much more quickly.
In particular, if the orbit is conjugate to a rigid rotation
\Eq{conjugacy} with a Diophantine rotation vector $\omega$
and the map $f$ and function $h$ are $C^\infty$, then \Eq{WB} converges faster than any power, \cite{Das18b}
\[
	|\WB_T(h) - \langle h \rangle|  < \frac{C_k}{T^k}, \quad \forall k \in \bN .
\]
Note that the function $\Psi$ is not the only $C^{\infty}$ weighting function with these excellent convergence properties. Recent work of \citet{Ruth24} uses reduced rank extrapolation to optimize the 
choice of weighting function.

We estimate the error of the WBA for a given function $h$ and a given time $T$ by
computing the effective number of digits of accuracy:
\beq{digits}
	\digT = -\log_{10} \left|\WB_{T}(h)(z)-\WB_{T}(h)(f^T(z)) \right| ,
\eeq
i.e., comparing the result for the first $T$ iterates with that for the next $T$ iterates.\cite{Sander20, Meiss21} 

To obtain a criterion distinguishing chaotic and regular orbits, 
we need to make a choice of cutoff value for $\digT$, declaring that orbits with
\beq{ChaosThreshold}
	\digT < D_T \Rightarrow \mbox{ ``chaotic'' } ;
\eeq
conversely, all orbits with $\digT \ge D_T$ are ``nonchaotic''.
Based on \Fig{ArnoldHist} for the 1D case, we choose the cutoff:
\beq{1DCriterion}
	 T = 10^5, \quad D_T = 9 \quad \mbox{(Circle Maps)}.
\eeq
This guarantees at 
least nine digits of accuracy in the computed rotation number. 
For the Arnold circle map in \Sec{arnold1D}, we know that there are no 
chaotic orbits for $0<a \le 1$; with the adopted criterion, our computations falsely identify only
$0.07\%$ ($1145$ out of $1.6$ million) of the orbits in this range to be chaotic.

For the 2D maps in \Sec{torusMaps}, we use the criterion
\beq{2DCriterion}
	T = 10^6, \quad D_T = 9 \quad \mbox{(Two-Torus Maps)}.
\eeq
Like the one-dimensional case, this is conservative in that chaotic orbits are quite unlikely to be identified as regular.

\subsection{Rational vs. Irrational: Farey Trees}\label{app:FareyTree}

In addition to providing the distinction between regular and chaotic orbits,
the WBA can be used to compute an accurate value of the time average of a function $h$.
In particular, we can compute the rotation vector \Eq{RotationVector} of an orbit for a torus 
map of the form \Eq{torusMap} using
 \beq{omegaWBA}
	\omega_T = \WB_T(F(x)-x) = \Omega +  \WB_T(g(x;a)).
\eeq
If $T$ is large enough and the rotation vector exists, we expect
$\omega_T \approx \omega(x,f)$ \Eq{RotationVector}. We note that this is by no means the only possible
useful $h$ to choose. For example, for quasiperiodic orbits, 
another choice of $h$ allow one to compute the conjugacy between
the map and a rigid rotation.\cite{Blessing23, Das16a}

In this section, we focus on circle maps, $d = 1$, and review the Farey
tree method of \citet{Sander20}, which when combined with \Eq{omegaWBA}
allows us to distinguish orbits that are periodic---those with rational
$\omega$---from those that are dense on a circle with irrational $\omega$. 

Even though a numerical determination of the irrationality of $\omega_T$ is impossible, we
will declare it to be ``effectively rational'' if it is sufficiently close to a \textit{low-order} $\tfrac{p}{q} \in \bQ$,
i.e., one with a ``small'' denominator $q$.
A rational is $\delta$-close to $x$ if it lies in the interval
\beq{Idelta}
	B_\delta(x) \equiv (x-\delta, x+\delta).
\eeq
The smallest denominator of a rational approximation within $\delta$ of $\omega$ is then
\beq{qmin}
	q_{min}(\omega,\delta) \equiv \min\{ q \in \bN : \tfrac{p}{q} \in B_\delta(\omega), p \in \bZ\} .
\eeq

We previously discussed an efficient method to compute \Eq{qmin}
that uses the Farey tree expansion of $\omega$.\cite{Sander20} The Farey (or Stern-Brocot) tree computes a sequence of 
rational approximations that converge to $\omega$.
We proved that the first such rational on the tree
that falls in $B_\delta(\omega)$ gives $q_{min}$.\cite{Sander20}
For example,
\[
	q_{min}\left( \sqrt{2},10^{-9} \right) = 33,461 
\]
(with corresponding numerator $p =  47,321$). Note that $q_{min} \sim 10^{4.52}$, 
which is close to $\frac{1}{\sqrt{\delta}} = 10^{4.5}$---as we see next, this is not unusual.

To decide what it means for a denominator to be ``small'', we computed
$q_{min}$ for a uniform distribution $\omega \in (0,1)$; since irrationals
have measure one, we take the computed distribution to be that of irrationals.
Numerically we found that the resulting log-denominators have a distribution
that is nearly symmetric about the mean
\beq{MeanDenom}
	\langle \log_{10} q_{min} \rangle = -\tfrac12 \log_{10}\delta + \alpha  ,
\eeq 
where $\alpha = - 0.05 \pm 0.001$, 
with standard deviation
\beq{StdDevDenom}
 \sigma = 0.2935 \pm 0.0006
\eeq
which is independent of $\delta$.  
In support of this, in a  recent paper \citet{Chen23} proved that as $\delta \to 0$,
$\log_{10} \langle q_{min} \rangle$ has the form \Eq{MeanDenom} with a different $\alpha$. 
However, this result is for the log of the mean, whereas our numerics were for the mean of the log. 
More recently, \citet{Marklof24b} has fully verified \Eq{MeanDenom}, showing that $\alpha = -0.0502959\dots$ and 
$\sigma = 0.293336\dots$, consistent with our numerical results.

It is important to note that it is not just small $q_{min}$ that correspond to a nearby rational: 
if an interval $B_\delta(\omega)$ is close to a low-order rational, but does not include this point, then $q_{min}$
can be much larger than the mean indicated by \Eq{MeanDenom}. For example, if 
$\omega$ is very close to $\tfrac{0}{1}$, but $\delta$ is so small 
that $0 \notin B_\delta(\omega)$, then the denominator can be unusually large; for example,
\[
	q_{min}(\sqrt{2}\times 10^{-8}, 10^{-9})  = 66,040,883  \sim 10^{7.8},
\]
which is $11\sigma$ above the mean \Eq{MeanDenom}.
Thus we will declare a number effectively irrational only if $q_{min}$ 
is close to the mean \Eq{MeanDenom} in the sense of the standard deviation $\sigma$.
Such an $\omega$ is ``typical'' in the sense of the uniform distribution and thus is ``irrational''.

For a given accuracy $\delta$, we declare that $\omega$ is an approximation of an irrational if
\beq{sTol}
	\left|\log_{10}(q_{min}(\omega,\delta)) + \tfrac12 \log_{10}(\delta) \right| < s, 
\eeq
for a given tolerance $s$. We typically choose
\beq{1DTolerance}
	\delta = 10^{-9}, \quad \mbox{and}  \quad s =  1.6875 \approx 5.75 \sigma.
\eeq
 This means that we declare that $\omega$ approximates an irrational for periods
\beq{IrrationalCriterion}
	649 < q_{min} < 1.54(10)^6 \Rightarrow \mbox{``irrational''}.
\eeq
This choice of tolerance $s$ means we are quite conservative in designating a rotation number as rational. 
Using Criterion \Eq{1DTolerance} we incorrectly identify $0.05\%$ of the random
$\omega$ as rational (since rationals have measure zero, the result should be zero).
Similarly, we see in \Sec{arnold1D} that \Eq{1DTolerance} 
erroneously identifies only $0.08\%$ of the nonchaotic orbits of the Arnold circle map as having irrational rotation number when $a \in (1,2.5)$, where it is known that there are no such orbits.

\subsection{Resonant vs. Incommensurate: Resonance Orders}\label{app:ResonanceOrder}
Here we recall a numerical method that generalizes the Farey tree method of \App{FareyTree} to higher dimensions.
In particular given $\omega_T$ we wish to compute the rank and resonance order, recall \Sec{ResOrdTheory}. 

A vector $\omega$ is approximately commensurate if $\left| m\cdot \omega -n \right|$ is 
small. In \citet{Meiss21} we developed a method for computing such commensurabilities. 
We say that a vector $\omega$ is $(m,n)$-resonant \textit{to precision $\delta$} if
the resonant plane \Eq{ResonantPlanes} intersects the ball \Eq{Idelta} about $\omega$, 
\beq{ResonancePlane}
	\cR_{m,n} \cap B_\delta(\omega) \neq \emptyset .
\eeq
Using the Euclidean norm, the minimum distance between the plane and the point $\omega$ is
\beq{minDist}
	\Delta_{m,n}(\omega) = \min_{\alpha \in \cR_{m.n}} \| \alpha - \omega\|_2 
	                     = \frac{|m \cdot \omega-n|}{\|m\|_2} .
\eeq
Thus $\omega$ is $(m,n)$ resonant to precision $\delta$, whenever $\Delta_{m,n}(\omega) < \delta$, and
we call the value
\beq{Mres}
	M(\omega,\delta) = \min \{\|m\|_1 : \Delta_{m,n}(\omega) < \delta, \, m \in \bZ^d\setminus\{0\}, \, n\in \bZ \},
\eeq
the \textit{resonance order} of $\omega$.

As far as we know, there is no generalization of the $d=1$ Farey tree result of \App{FareyTree}
to compute \Eq{Mres} efficiently.\footnote
{One could use the Kim-Ostlund tree to get resonance relations;\cite{Ashwin93} however, it is not clear that this algorithm returns a minimal $\|m\|$.} 
Nevertheless, since there are finitely many $m \in \bZ^d$ such that $\|m\|_1 \le M$, a brute force computation
is possible for modest values of $M$. \cite{Meiss21}

To understand what resonance orders are ``typical,'' we computed the minimal resonance order
\Eq{Mres} for a set of equi-distributed, random $\omega \in [0,1)^2$ as a function of the precision $\delta$.
\cite{Meiss21} The resulting distribution of $\log(M)$ has a mean\cite{Meiss21}
\beq{MeanRes}
	\langle \log_{10} M(\omega,\delta) \rangle = -0.334\log_{10}(\delta) -0.091.
\eeq
\citet{Marklof24a, Marklof24b} has recently obtained theoretical results that are consistent with \Eq{MeanRes}.
We find a  standard deviation of
\[
	\sigma = 0.171.
\]
As in the one-dimensional case, the standard deviation seems to be essentially independent of $\delta$.

Since the cutoff \Eq{2DCriterion} gives rotation number calculations 
accurate to within $10^{-9}$, we choose $\delta = 10^{-9}$. For this case, \Eq{MeanRes} implies 
that $\langle\log_{10} M\rangle =  2.915$. We declare that a vector is nonresonant if 
\beq{Incommensurate}
	256 \le M \le 2673 \Rightarrow \mbox{``nonresonant''},
\eeq
corresponding to $ 2.407 < \log_{10}(M) < 3.427 $,
which is a range of approximately $\pm 3 \sigma$ about the mean \Eq{MeanRes}.
To test this criterion we selected $10^4$ randomly distributed values 
uniformly in $[0,1)^2$, and found that $1.36\%$ were incorrectly identified as resonant. 
Note that the distribution of log-orders for random vectors is not symmetric around the mean; 
in particular, $M<256$ occurred $1.32\%$ of the time, and $M>2673$ occurred $0.04\%$ of the time. 

We can further categorize the orbits that are determined to be resonant (those 
that fail Criterion \Eq{Incommensurate}) by the rank of the resonance. 
Rank-two resonant orbits have frequencies on the intersection of a pair of different resonance lines,
recall \Fig{ResLines}. That is, both of the components of the rotation vector are ``rational''.
These can be identified using the criterion
\beq{Periodic}
	M \mbox{ fails  \Eq{Incommensurate} and } 
	\omega_1, \omega_2 \mbox{ fail \Eq{IrrationalCriterion} } 
	\Rightarrow \mbox{``periodic''}.
\eeq

If a resonant orbit is not periodic, then it lies on a single resonance line and so has rank one.
We will simply refer to such orbits as ``resonant''; they typically are dense on topological circles.
Thus the criterion for a (rank-one) resonant orbit is
\beq{Resonant}
	M \mbox{ fails \Eq{Incommensurate} but at least one of } \omega_1,\omega_2 \mbox{ satisfy \Eq{IrrationalCriterion}}
	 \Rightarrow \mbox{``resonant''}.
\eeq
The criteria \Eq{Periodic} and \Eq{Resonant} are used to distinguish the orbit types for the 2D maps in \Sec{torusMaps}.

\section{Critical $\eps$}\label{app:EpsCrit}

In this appendix, we explain why $\epsC$ \Eq{epsCrit} almost always
exists for most choices of amplitudes $a_i$ and phases $\phi_i$. 
We first argue that typically there is some point $(x_1,x_2)$ for which $\det(H) < 0$ for
the matrix $H$ in \Eq{defineH}. First, if we choose $x_2 = -\phi_2$, then 
\[
	\det(H(x_1,-\phi_2)) = a_1 a_4 \sin(2 \pi(x_1+\phi_1))\sin(2 \pi(\phi_4-\phi_2)).
\]
Therefore if $a_1 a_4 \neq 0$ and $\phi_4-\phi_2 \neq n \pi$ for some integer $n$,
this determinant is nonzero and odd about $x_1 = -\phi_1$.
Thus $\det(H) < 0$ at some point.
If $\phi_4-\phi_2 = n \pi$, we can consider a similar argument upon choosing $x_1 = -\phi_3$ as long as $\phi_3 - \phi_1 \neq m \pi$ for some integer $m$. Equivalent arguments apply if $a_2a_3 \neq 0$, as long as $\phi_4-\phi_2 \neq n \pi$ or 
$\phi_3 - \phi_1 \neq m \pi$. 

An exceptional case would be if both $\phi_4-\phi_2 = n \pi$, $\phi_3 - \phi_1 = m \pi$.
Now choose $x_1 =  \tfrac{\pi}{2}-\phi_1 $ and $x_2 =\tfrac{\pi}{2}-\phi_4 $. Then 
\bsplit{SpecialDetH}
	\det H &= a_1 a_4 \sin \tfrac{\pi}{2} \sin \tfrac{\pi}{2} - a_2 a_3 \sin \left( \tfrac{\pi}{2} + 
		\phi_3-\phi_1 \right) \sin \left(\tfrac{\pi}{2} + \phi_2 - \phi_4 \right)\\
		&= a_1 a_4 - a_2 a_3 \sin\left(\tfrac{\pi}{2} + n \pi \right) \sin \left(\tfrac{\pi}{2} - m \pi \right)
\esplit
Note that if $m$ and $n$ are both even or both odd (i.e., have the same parity), then at this point $H = a_1 a_4 - a_2 a_3$. 
Thus long as $a_1 a_4 - a_2 a_3 < 0$, $\det(H) < 0$ at this point. On the other hand if $a_1 a_4 - a_2 a_3 > 0$, then use $x_2 = \tfrac{3 \pi}{2}-\phi_2$, which flips the signs of both terms, again giving $\det(H) < 0$.
If however, $m$ and $n$ have the opposite parity, then 
at the point \Eq{SpecialDetH} $\det(H) = a_1 a_4 + a_2 a_3$. If this is negative, we are done. Otherwise
use $x_2 = -\phi_2 + \tfrac{3\pi}{2}$, which flips the signs of both terms.

If none of the above cases hold, then $a_1 a_4 - a_2 a_3 = a_1 a_4 + a_2 a_3 = 0$, $\phi_4-\phi_2 = n \pi$, $\phi_3 - \phi_1 = m \pi$, implying
$\det(H)$ is identically zero. But as long as the trace is nonzero, we can use a similar argument for the term that is linear in $\eps$
in \Eq{defineH} to show that $\epsC$ exists.

Therefore the only exception to the existence of $\epsC$ is the case $a_1=a_4 = 0$ and $a_2a_3 = 0$,
which gives trivial dynamics. 

\newpage
\section{Parameters for \Sec{torusMaps}}\label{app:parameters}

This appendix gives the parameters for the computations in \Sec{torusMaps}. 
The parameters $\eps$ and $\Omega$ for the images in \Fig{orbitmontage} are shown in \Tbl{orbitmparamsB}. 
This also gives the computed $\omega_T$ for the regular orbits---cases (a)-(d),
and the precision $\digT$ for each case.
\Tbl{orbitmparamsC} gives amplitudes and phases for \Eq{forceFunction}
used in \Figs{orbitmontage}{singleOmega} (\Case{0}) and in \Fig{manycases} (eight cases)
along with the calculated values of $\epsC$. 
For all of these parameter sets, $\|a\|_1 = 1$.

\begin{table}[ht]
\begin{center}
\begin{tabular}{c l c c r}
\hline
Label & $\eps$   & $\Omega$ & $\omega_T$ & \multicolumn{1}{c}{$\digT$} \\
\hline   
   (a) & 0.8   & (0.2, 0.7)  &$(0.19570941533509,0.70456941766774)$ &$14.7505$\\
  (b) & 0.8  & (0.84, 0.835) &$(0.83947029089469, 0.83947029089470)$ &$15.2556$\\
  (c) & 0.8  & (0.5,   0.7)  &$(0.49778852806059, 0.70330086015610)$ &$12.0529$\\
  (d) & 1.5  & (0.1,   0.8)  &$(0.07425024047212, 0.83592850272224)$ &$14.3471$\\
  (e) & 2.6  & (0.7,   0.3) & &$4.0484$\\
  (f) & 4.0  & (0.24,   0.4)& & $2.5522$\\
\hline\hline
\end{tabular}
\end{center}
\caption{The parameters $\eps$ and $\Omega$ for the orbits in 
\Fig{orbitmontage}, using the force \Eq{forceFunction}
with amplitudes and phases for \Case{0} in \Tbl{orbitmparamsC}.
The last two columns give the computed $\omega_T$ and $\digT$ using $T=10^6$.} 
\label{tbl:orbitmparamsB}
\end{table}

\begin{table}[ht]
\begin{center}
\begin{tabular}{c c l l l l l}
\hline
Case & Param. & $1$ & $2$ & $3$ & $4$ & $\epsC$ \\
\hline
\multirow{2}*0 & $a_i$	& $0.221320306832860$  & $0.220593736048273$ &   $0.152270586812051$ & $0.405815370306816$ & \multirow{2}*{$2.22044$} \\
&  $\phi_i$ & $0.369246781120215$  &  $0.111202755293787$  & 	$0.780252068321138$ & 	$0.389738836961253$\\
\hline 
\multirow{2}*1 & $a_i$  & $ 0.406588842221655$ & $ 0.062715680327705$ &  $0.179066359898821$ & $0.351629117551819$ & \multirow{2}*{$2.2070$}\\
 & $\phi_i$ & $0.957506835434298$ &  $ 0.964888535199277$ &  $  0.157613081677548$ & $0.970592781760616$\\
\hline
\multirow{2}*2 & $a_i$ & $0.211681398612178 $& $  0.317651811580494$ &  $ 0.375591536887180 $&   $0.095075252920149$ &  \multirow{2}*{$2.4566$}\\
&$\phi_i$ & $0.273022072458714$ &  $ 0.542430207288253$ &  $ 0.431224181579691 $& $0.153093675447227$ \\
\hline
\multirow{2}*3 & $a_i$ & $0.012536281513538$ &    $0.465737538631897$ &  $  0.503609970119032$ &    $0.018116209735533$ & \multirow{2}*{$2.0564$}\\
&$\phi_i$  & $0.739790415703666$ &    $0.023926884448995$ &    $0.490328482174893$ &    $0.304888898615625$\\
\hline
\multirow{2}*4 & $a_i$  & $0.760566444256527 $&               0 &                  0 &  $ 0.239433555743473 $ &  \multirow{2}*{$1.3148$} \\
&$\phi_i$ & $0.739790415703666$ &  $0.023926884448995$ &   $0.490328482174893$ &   $0.304888898615625$ \\
\hline
\multirow{2}*5 & $a_i$ &              0 &   $0.760566444256527$ &  $ 0.239433555743473$ &                  0 &  \multirow{2}*{$2.3434$} \\
&$\phi_i$ &$ 0.739790415703666$ &  $ 0.02392688444899$ &  $ 0.490328482174893$ &   $0.304888898615625$ \\
\hline
\multirow{2}*6 & $a_i$ & $0.007280035519179 $&   $0.942703650246408$ &   $0.039647117954398$ &   $0.010369196280015$  & \multirow{2}*{$5.2100$}\\
&$\phi_i$ & $0.384398913909761$ &   $0.203897175276146$ &   $0.913862879483257$ &   $0.191420654770675$ \\
\hline
\multirow{2}*7 & $a_i$ &                0 &   $0.352156017226267$ &                  0 &   $0.647843982773733 $ &  \multirow{2}*{$1.5436$}\\ 
&$\phi_i$ & $0.369246781120215$  &   $0.111202755293787$  &  $0.780252068321138$ &   $0.389738836961253$ \\
\hline\hline
\end{tabular}
\end{center}
\caption{Amplitudes, phases, and $\epsC$ for curves in \Figs{orbitmontage}{singleOmega} (\Case{0}) and in \Fig{manycases}
(\textit{Cases} (0-7)). }
\label{tbl:orbitmparamsC}
\end{table}

\newpage
\bibliography{TorusMaps}

\end{document}